\documentclass[12pt]{article} 

\usepackage[left= 0.8in, textwidth=6.9in, textheight=9.0in, top=1.0in]{geometry}
\usepackage{graphicx}
\usepackage{framed}
\usepackage{subcaption}
\usepackage{amsmath,amsfonts,amssymb,verbatim}
\usepackage{enumitem}
\usepackage{color}  
\usepackage{bm}                                  

\newcommand{\nn}{\nonumber}
\newcommand{\MT}{\left[ \begin{array}{rrrrrrrrrrrrrrrrrrrr}}
\newcommand{\EM}{\end{array}\right]}
\newcommand{\EQ}{\begin{equation}\begin{alignedat}{4}} 
\newcommand{\EE}{\end{alignedat}\end{equation}}
\newcommand{\SEQ}{\begin{subequations}}
\newcommand{\ESE}{\end{subequations}} 

\newcommand{\Real}{\mathbb R}
\newcommand\norm[1]{\left\lVert#1\right\rVert}

\newtheorem{theorem}{Theorem}
\newtheorem{lemma}{Lemma}

\newtheorem{definition}{Definition}

\newtheorem{corollary}{Corollary}

\def\Fr{\ds \frac}
\def\ds{\displaystyle}

\def\calI{{\cal I}}

\def\bff{{\bf f}}
\def\bfe{{\bf e}}

\def\bfh{{\bf h}}

\def\bfu{{\bf u}}

\def\bfy{{\bf y}}

\title{\LARGE \bf
A Surrogate Data Assimilation Model for the Estimation of Dynamical System in a Limited Area \thanks{This work was supported in part by U.S. Naval Research Laboratory - Monterey, CA and NSF \# 2202668.}
}

\author{Wei Kang\thanks{Department of Applied Mathematics, Naval Postgraduate School, Monterey, CA, USA; wkang@nps.edu, hzhou@nps.edu}  \thanks{Department of Applied Mathematics, University of California, Santa Cruz, CA, USA}
\and Liang Xu\thanks{Marine Meteorology Division, Naval Research Laboratory, Monterey, CA, USA; liang.xu@nrlmry.navy.mil}
\and Hong Zhou$^\dagger$
}

\begin{document}
\maketitle

\abstract{We propose a novel learning-based surrogate data assimilation (DA) model for efficient state estimation in a limited area. Our model employs a feedforward neural network for online computation, eliminating the need for integrating high-dimensional limited-area models. This approach offers significant computational advantages over traditional DA algorithms. Furthermore, our method avoids the requirement of lateral boundary conditions for the limited-area model in both online and offline computations. The design of our surrogate DA model is built upon a robust theoretical framework that leverages two fundamental concepts: observability and effective region. The concept of observability enables us to quantitatively determine the optimal amount of observation data necessary for accurate DA. Meanwhile, the concept of effective region substantially reduces the computational burden associated with computing observability and generating training data.} 

\section{Introduction}
Data assimilation (DA) is a highly interdisciplinary field that merges estimation theories and stochastic systems with scientific computing algorithms. Its primary objective is to integrate observations with numerical models of dynamical systems, aiming to achieve an optimal estimation of the system's state or parameters as it evolves over time. The applications of DA span across a diverse range of scientific and engineering domains. Notable examples include numerical weather prediction, ocean forecasting, space weather analysis, robotics, and combustion dynamics, among many others. Various DA algorithms have been extensively explored in the literature to cater to a wide range of applications. Notably, two families of prominent methods, namely ensemble Kalman filters (EnKF) and 4D-variational methods (4DVar), have demonstrated remarkable success in the realm of numerical weather prediction. For a comprehensive understanding of these algorithms, interested readers are encouraged to delve into the works of \cite{houtekamerzhang,bannist1,bannist2,fairbairnpring2014,vetracarvalhovanleeuwen2018,park_zupanski} on EnKF and \cite{park_zupanski,rabier2005,lorencbowler2015,xurosmond2005,fairbairnpring2014} on 4DVar and the references therein. It is important to note that DA plays a pivotal role in diverse engineering applications as well, for example robotic systems \cite{changedwards2019}, combustion dynamics \cite{yujaravel2019}, shock wave analysis \cite{srivastavakang2023}, and control systems \cite{lagordavis2016}. In recent years, learning-based methods for nonlinear filtering have garnered considerable attention. Specifically, the findings presented in this paper exhibit a similar structure to the deep filters introduced in \cite{qianyin2023, kangxu2022}, with the exception that our focus lies in the limited-area estimation of systems defined by partial differential equations (PDEs).

Dynamical phenomena with highly variable nature, such as water vapor in cloud and precipitation formation, are often inadequately represented in global models. To overcome this limitation, DA techniques employing mesoscale or microscale limited-area models (LAMs) have been utilized, enabling the acquisition of high-resolution results that would be computationally prohibitive with a global model \cite{davies2014}. Several algorithms, including EnKF and 4DVar, have been extensively studied for limited-area DA \cite{mengzhang2011, dowellzhang2004, tornhakim2008, gustafssonhuang2012}. While each approach has its advantages and disadvantages, existing algorithms for limited-area DA face common challenges. Firstly, the running of LAMs requires  lateral  boundary  conditions, which are  typically  obtained  from a global model or another lower-resolution LAM containing the area of interest. This requirement not only increases the computational burden but also introduces potential errors into the LAM due to the lateral boundary conditions \cite{gustafssonkallen1998}. Secondly, current DA algorithms necessitate multiple iterations of running the numerical model within each estimation cycle. This computational intensity is the primary reason why DA for high dimensional models is computationally expensive.

In this paper, we introduce a learning-based surrogate data assimilation model, referred to as the surrogate DA model, for state estimation in a limited area. Our model employs deep learning techniques, where a neural network is trained using sensor data as inputs and providing the estimated values of the system's states as outputs. Our approach offers numerous benefits. Firstly, computationally intensive tasks, such as generating data and training the neural network, are all performed offline during the design phase of the surrogate DA model. This leads to a significant reduction in online computational burden as the evaluation of the surrogate DA model does not require propagating the dynamic model. This is fundamentally different from conventional DA algorithms, such as EnKF or 4D-Var, that require propagating the dynamic model multiple times in every estimation cycle.  Secondly, our method does not rely on lateral boundary conditions around the limited area in both online and offline computations. Moreover, the design of the surrogate DA model is based on a solid theoretical foundation that employs two fundamental concepts: observability and effective region. The concept of observability is used to quantitatively answer a key question of learning: how much observational data is enough for the purpose of DA? The concept of effective region, on the other hand, greatly reduces the computational load involved in computing observability and generating training data.

The concept of observability and its computational methods are introduced in Section \ref{sec_obs} for general discrete-time dynamical systems. Section \ref{sec_obs_effective} then extends this concept to discretized PDEs on a grid, introducing the notion of an effective region and providing illustrative examples. Next, in Section \ref{sec_learning}, we present a learning-based approach for the design of surrogate DA models. The effectiveness of this approach is demonstrated using the Burgers equation in $\mathbb{R}^2$ as a case study. Furthermore, in Section \ref{sec_prediction}, we apply the learning-based method to develop a surrogate prediction model. It is illustrated with an example of tsunami forecasting. 

\section{Observability}
\label{sec_obs}
Consider a discrete-time dynamical system
\SEQ
\label{eq_sys}
\begin{alignat}{4}
&\bfu (k+1)=\bff (\bfu(k)) , &\quad& \bfu \in \Real^n, &\quad& k\in {\mathbb N}_0, \label{eq_sysa}\\
&\bfy(k)=\bfh (\bfu (k))+\boldsymbol\nu(k), && \bfy \in \Real^m,  \label{eq_sysb}\\
&z(k)=w(\bfu (k)), && z \in \Real,
\end{alignat}
\ESE
where $\bff: \Real^n\rightarrow \Real^n$, $\bfh: \Real^n\rightarrow \Real^m$  and $w: \Real^n\rightarrow \Real$ are Lipschitz continuous functions, $\bfu$ is the state variable of the system, $\bfy$ is the output, or observation, whose value is known or can be measured, $\boldsymbol\nu$ is the observation uncertainty, $z$ is the variable to be estimated. In this section, we develop a theoretical framework to answer the following question.\\
\vspace{-0.125in}

\noindent
\textbf{Question}. Let $[0, K]$ be a time interval.  Given $\{ \bfy(k)\}_{k=0}^K $, is it possible to accurately approximate the value of $z(K)$? \\
\vspace{-0.125in}

This problem is closely linked to the concept of observability. In control theory, observability is defined based on the rank condition of an observability matrix. However, this definition fails to encompass some crucial aspects pertinent to high-dimensional systems. In the following, we present a quantitative measure of observability for individual state variables. Unlike the classical control theory definition, where observability is defined for the entire system, the introduced quantitative measure in this section allows for separate treatment of different state variables. This becomes particularly important when computing the observability for the entire system is unnecessary in limited area DA. Moreover, this quantitative measure can be numerically computed based on a matrix with reduced dimensions, while in control theory, the observability matrix has the same dimension as the system. Most importantly, the concept and findings presented in this section lead to a novel concept known as the ``effective region," which is crucial to the learning-based surrogate DA model that will be introduced later in this paper.

\subsection{Linear systems}
Consider a linear system,
\SEQ
\label{eq_linsys}
\begin{alignat}{4}
&\bfu(k+1)=A\bfu(k), &\quad& \bfu \in \Real^n, &\quad& k\in {\mathbb N}_0,\label{eq_linsysa}\\
&\bfy(k)=H\bfu(k)+\boldsymbol\nu(k), && \bfy \in \Real^m,\\
&z(k)=W\bfu (k), && z \in \Real, \label{eq_linsysc}
\end{alignat}
\ESE
where $A\in\Real^{n\times n}$, $W\in\Real^{1\times n}$, $H\neq 0$ and $H \in \Real^{m\times n}$. In control theory, if the state $\bfu(K)$ at any time $K$ can be uniquely determined by $\{ \bfy(k)\}_{k=0}^K$, then the system is said to be {\it observable}. In this paper, we expand the meaning of this concept by introducing a broader perspective. Specifically, if $\{ \bfy(k)\}_{k=0}^K$ contains adequate information to accurately determine $z(K)$, we say that $z(K)$ is observable. This concept is mathematically defined later in this section.  

When calculating observability, we make the assumption that $z$ is a scalar. There are multiple reasons for computing observability on a per-variable basis. In systems with high dimensions, the observability of state variables can vary greatly: some variables are strongly observable, some are weakly observable, and some may even be unobservable. Furthermore, by analyzing each variable independently, we can concentrate on the observability of those state variables that are relevant to a specific area of interest, without having to consider the overall state variable. This approach results in significant dimension reduction in computation. Moreover, for a discretized PDE system, the observability of the state variables at nearby grid points is close to each other. Therefore, the distribution of the PDE's observability in a region in space can be approximated by computing the observability at a set of sample points. Furthermore, it will be proved in this section that the quantity that measures the observability of scalar variables has a compact formula. 

Given $\bfy(0),\bfy(1),\cdots,\bfy(K)$, where $K\geq 1$ is an integer, the goal is to estimate the scalar $z(K)$ in (\ref{eq_linsysc}).
For example, if $W=[ 1, 0, \cdots, 0]$, the variable to be estimated is $u_1(K)$. Consider two solutions of (\ref{eq_linsys}): $\{ \bfu (k)\}_{k=0}^K$ and $\{\hat \bfu (k)\}_{k=0}^K$.
The values of $\{ \bfy(k)\}_{k=0}^K$ and $\{ \hat \bfy(k)\}_{k=0}^K$ may or may not contain enough information to determine whether the solutions are different. This is especially true in the presence of sensor errors. If $\bfy(k)-\hat \bfy(k)$ is within the magnitude of sensor error, an estimator cannot determine which solution is the true one. In the following, we introduce a concept about distinguishability that is closely related to the observability of the system's solutions.

\begin{definition}
Let $\epsilon >0$ be a real number. Two solutions, $\{ \bfu (k)\}_{k=0}^K$ and $\{\hat \bfu (k)\}_{k=0}^K$, of (\ref{eq_linsys}) are said to be $\epsilon$-indistinguishable, or simply indistinguishable,  if 
\EQ
\label{eq_indist}
\ds\sum_{k=0}^K\norm{H\hat\bfu(k)-H\bfu(k)}_2^2 \leq \epsilon^2.
\EE
Any $W\hat \bfu (K)$ in which $\hat  \bfu (K)$ satisfies (\ref{eq_indist}) is called an indistinguishable estimate of $z(K)$. 
\end{definition}
In this definition,  $\epsilon>0$ is a constant that represents an upper bound of observation uncertainty. If the set of indistinguishable estimates of $z(K)$ has a small diameter, then any estimate from this set has a small error. Motivated by the work on observability in \cite{kangxu2009}, we define a constant $\rho$ that is a quantitative measure of this diameter. Given a nominal solution satisfying (\ref{eq_linsys}) and $\epsilon >0$, we define
\begin{subequations}
\label{def1}
\begin{alignat}{4}
&\rho^2=\ds\max_{\hat\bfu(0),\cdots,\hat\bfu(K)} \{ ( W \hat \bfu(K)-z(K))^2\}, \label{def1a}\\ 
&\hspace{-0.2in} \mbox{subject to} \nn \\
&\hat \bfu(k+1)=A\hat \bfu(k),\;\;\;\;\;\;\; k=0,1,\cdots,K-1, \label{def1b}\\
&\ds\sum_{k=0}^K\norm{H\hat\bfu(k)-H\bfu(k)}_2^2 \leq \epsilon^2.  \label{def1e}
\end{alignat}
\end{subequations}

The value of $\rho$ represents the maximum potential error in estimating $z(K)$ among all indistinguishable estimates. When using estimators that leverage additional information beyond observation data, such as the Kalman filter that incorporates both observations and error covariance, they have the potential to produce estimates that surpass the worst-case estimation from the indistinguishable set. Therefore, the error of a good estimator should be smaller than $\rho$. For a given error tolerance, $tol > 0$, all indistinguishable estimates are within the $tol$-neighborhood around $z(K)$ if $\rho < tol$. In this case, $z(K)$ is said to be  strongly observable.

To compute $\rho$, we reformulate (\ref{def1}) as follows. A solution can be uniquely determined by its initial state $\hat\bfu_0$,
\EQ
\label{eq_traj1}
&\hat \bfu(k)= A^k \hat\bfu_0,\\
& H\hat \bfu(k)=HA^k\hat \bfu_0.
\EE
The following matrix is essential in the computation of $\rho$,
\EQ
\label{eq_Gramian}
G=\ds\sum_{k=0}^K (A^\intercal)^kH^\intercal HA^k.
\EE
In control theory, $G$ is called the observability Gramian. It is a symmetric matrix that measures the observability of control systems \cite{kangxu2009,kailath}. We denote the initial error by 
\EQ
\label{eq_Du}
\varDelta \bfu_0 = \hat \bfu_0-\bfu_0.
\EE
Substituting (\ref{eq_traj1}) into (\ref{def1}) yields a quadratically constrained quadratic programming (QCQP). 
\begin{lemma}
\label{lemma1}
The constrained maximization in (\ref{def1}) is equivalent to the following QCQP
\SEQ
\label{eq_quadratic}
\begin{alignat}{4}
&\rho^2=\ds\max_{\varDelta \bfu_0} \varDelta \bfu_0^\intercal (W(K))^\intercal W(K) \varDelta \bfu_0,\\
&\hspace{-0.2in} \mbox{subject to} \nn\\
&\varDelta \bfu_0^\intercal G \varDelta \bfu_0=\epsilon^2, \label{eq_quadraticb}
\end{alignat}
\ESE
where $W(K)=WA^K$ and $G$ is the observability Gramian (\ref{eq_Gramian}).
\end{lemma}
The proof of this lemma is straightforward. Substituting (\ref{eq_traj1}), (\ref{eq_Gramian}), (\ref{eq_Du}) into (\ref{def1}) yields (\ref{eq_quadratic}). The inequality in (\ref{def1e}) is replaced by an equation, (\ref{eq_quadraticb}), because the maximum value is always achieved on the boundary. Based on this lemma, the computation of  $\rho$ boils down to solving (\ref{eq_quadratic}). For problems that have high dimensions, finding a numerical solution to (\ref{eq_quadratic}) is not as straightforward as it may seem if the overall system is weakly observable or unobservable. In this case, the condition number of $G$ is very large. In the following, we first prove a theorem in which $G$ has full rank. The case in which $G$ is singular will be addressed in Section \ref{sec_singularG}. 

\begin{theorem}
\label{thm1}
Let $\sigma_1, \sigma_2, \cdots, \sigma_n$ be the eigenvalues of $G$, the observability Gramian defined in (\ref{eq_Gramian}). Let $T\in \Real^{n\times n}$ be a matrix in which the $i$th column is a unit eigenvector associated with $\sigma_i$, $1\leq i\leq n$. Suppose that $G$ has full rank. Then 
\EQ
\label{eq_unobs_formular}
\rho^2 = \epsilon^2 \ds\sum_{i=1}^n{\Fr{\bar w_i^2}{\sigma_i}},
\EE
where $\bar w_i$ is the $i$th component of the vector
\EQ
\label{eq_sbar}
\bar W=W(K) T.
\EE
\end{theorem}

{\it Proof}. To solve (\ref{eq_quadratic}), define the Lagrangian
\EQ
L(\varDelta \bfu_0,\lambda)=\varDelta \bfu_0^\intercal (W(K))^\intercal W(K) \varDelta \bfu_0-\lambda (\varDelta \bfu_0^\intercal G \varDelta \bfu_0-\epsilon^2),\nn
\EE
where $\lambda\in \Real$ is the Lagrange multiplier. At the solution of  (\ref{eq_quadratic}), denoted by $\varDelta \bfu_0^\ast$, there exists $\lambda^\ast$ such that
\EQ
&\Fr{\partial L}{\partial \varDelta \bfu_0} (\varDelta \bfu_0^\ast,\lambda^\ast)=0. \nn
\EE
Therefore, we have
\SEQ
\begin{alignat}{4}
&(W(K) )^\intercal W(K)  \varDelta \bfu_0^\ast=\lambda^\ast  G \varDelta \bfu_0^\ast, \label{eq_exp1}\\
&(\varDelta \bfu_0^{\ast})^\intercal G \varDelta \bfu_0^\ast=\epsilon^2.
\end{alignat}
\ESE
Multiplying $(\varDelta \bfu_0^\ast)^\intercal$ to (\ref{eq_exp1}), we have
\EQ
\rho^2&=(\varDelta \bfu_0^\ast)^\intercal (W(K) )^\intercal W(K)  \varDelta \bfu_0^\ast\nonumber \\
&=\lambda^\ast  (\varDelta \bfu_0^\ast)^\intercal G \varDelta \bfu_0^\ast\nn\\
&=\lambda^\ast \epsilon^2 .
\EE
Therefore, 
\EQ
\label{eq_lmda}
&\lambda^\ast = \rho^2/\epsilon^2.
\EE
Because $T$ is the matrix of unit eigenvectors and $G$ is symmetric, we have
\EQ
\label{eq_eigG}
 T^\intercal GT=diag\left(\MT \sigma_1 & \sigma_2&\cdots &\sigma_n \EM\right). \nonumber
\EE
Multiplying $T^\intercal$ to (\ref{eq_exp1}), applying $\bar W$ in (\ref{eq_sbar}) and the fact that $T^\intercal$ is the inverse matrix of $T$, yield
\EQ
\bar W^\intercal \bar W T^\intercal \varDelta u_0^\ast = \lambda^\ast diag\left(\MT \sigma_1 & \sigma_2&\cdots &\sigma_n \EM\right)T^\intercal\varDelta u_0^\ast. \nonumber
\EE
Because $ \bar W T^\intercal \varDelta u_0^\ast$ is a scalar, we have
\EQ
\lambda^\ast T^\intercal \varDelta u_0^\ast = (\bar W T^\intercal \varDelta u_0^\ast)\left[\begin{array}{cccc} \Fr{\bar w_1}{\sigma_1}& \cdots & \Fr{\bar w_n}{\sigma_n}\end{array}\right]^\intercal. \nonumber
\EE
Multiplying this equation by $\bar W$, we have
\EQ
\label{eq_exp2}
\lambda^\ast \bar W T^\intercal \varDelta u_0^\ast = (\bar WT^\intercal \varDelta u_0^\ast)\bar W\left[\begin{array}{cccc} \Fr{\bar w_1}{\sigma_1}& \cdots & \Fr{\bar w_n}{\sigma_n}\end{array}\right]^\intercal .
\EE
If $\bar W T^\intercal \varDelta u_0^\ast\neq 0$, then (\ref{eq_unobs_formular}) is proved by
cancelling the nonzero scalar  $\bar W T^\intercal \varDelta u_0^\ast$ in (\ref{eq_exp2}) and applying (\ref{eq_lmda}). If $\bar W T^\intercal \varDelta u_0^\ast = 0$, then
\EQ
0&=(\bar W T^\intercal \varDelta u_0^\ast)^2\nonumber\\
&=(W(K) \varDelta u_0^\ast)^2\\
&=\rho^2.
\EE
Because $\rho=0$, (\ref{eq_quadratic}) implies $W(K)=0$ and $\bar W=0$. Therefore, (\ref{eq_unobs_formular}) holds true.
$\blacksquare$

\subsection{Singularities in the observability Gramian}
\label{sec_singularG}
In applications, the eigenvalues of $G$ could be zero or extremely small. This is especially true for high dimensional systems that  have only a small number of observations. To remedy this problem, we modify $G$ by adjusting its eigenvalues based on the upper bound of initial state error, and the resulting matrix has full rank. Let $\sigma_1, \sigma_2, \cdots, \sigma_n$ be the eigenvalues of $G$ and let $T\in \Real^{n\times n}$ be the matrix formed by the corresponding unit eigenvectors. We have
\EQ
G=T \mbox{diag} (\sigma_1, \sigma_2, \cdots, \sigma_n ) T^\intercal .
\EE
If there is a zero eigenvalue, i.e., $\sigma_i=0$ for some $i$, then  for any $\varDelta \bfu_0$ that is a scalar multiple of $T_i$, the $i$th column vector in $T$, we have
\EQ
\label{eq_duTi}
\varDelta \bfu_0^\intercal G \varDelta \bfu_0 =0 < \epsilon^2 .
\EE
According to Lemma \ref{lemma1}, (\ref{eq_duTi}) implies that $u_0+\varDelta \bfu_0$ is indistinguishable from $\bfu_0$, no matter how large $\norm{\varDelta u_0}_2$ is. In reality, however, the error of initial guess cannot be arbitrarily large. If we assume that the unknown initial state is within a $\delta $-neighborhood of $\bfu_0$ for some $\delta > 0$, then in the direction of $T_i$, $\bfu_0+\varDelta\bfu_0$ is considered distinguishable from $\bfu_0$ if $\norm{\varDelta \bfu_0}_2 > \delta$. To reflect the known radius of the set of initial states, we introduce the following modified Gramian, 
\EQ
\label{eq_Gmodified}
&G_\delta = T \mbox{diag}(\tilde \sigma_1, \tilde \sigma_2, \cdots, \tilde \sigma_n) T^\intercal ,\\
&\tilde \sigma_i = \max\left\{\sigma_i, \epsilon^2/\delta^2 \right\}, & i=1, 2, \cdots, n.
\EE
If  $\varDelta \bfu_0$ is a scalar multiple of $T_i$, then $\bfu_0+\varDelta\bfu_0$ is  indistinguishable from $\bfu_0$, i.e.,
\EQ
  \varDelta \bfu_0^\intercal G_\delta \varDelta\bfu_0 = \tilde \sigma_i \varDelta\bfu_0^\intercal \varDelta\bfu_0< \epsilon^2 , \nonumber
\EE
implies that $\norm{\varDelta\bfu_0} < \delta$.  It is worth noting that $G_\delta$ does not have eigenvalues that are smaller than $\epsilon^2/\delta^2$. Replacing $G$ in (\ref{eq_quadraticb}) by the modified observability Gramian yields the following nonlinear programming, referred to as the primary QCQP, 
\SEQ
\label{eq_quadratic1}
\begin{alignat}{4}
&\rho^2=\ds\max_{\varDelta \bfu_0} \varDelta \bfu_0^\intercal (W(K))^\intercal W(K) \varDelta \bfu_0,\\
&\hspace{-0.2in} \mbox{subject to} \nn\\
&\varDelta \bfu_0^\intercal  G_\delta \varDelta \bfu_0=\epsilon^2. \label{eq_quadratic1b}
\end{alignat}
\ESE

\noindent Applying Theorem \ref{thm1}, the solution of the primary QCQP has an explicit solution. 
\begin{corollary}
\label{corollary1}
Let $G$ and $G_\delta$ be the observability Gramian (\ref{eq_Gramian}) and the modified observability Gramian (\ref{eq_Gmodified}), respectively. The eigenvalues of $G$ are $\sigma_1, \sigma_2, \cdots, \sigma_n$. Their corresponding unit eigenvectors form a matrix, $T$. Let $\rho$ be the maximum value of the primary QCQP (\ref{eq_quadratic1}). 
Then
\EQ
\label{eq_tilderho}
 \rho^2 =  \ds\sum_{i=1}^n\bar w_i^2 \min\{\Fr{\epsilon^2}{\sigma_i}, \delta^2\},
\EE
where $\bar w_i$ is the $i$th component of the vector
\EQ
\bar W=W(K) T.
\EE
\end{corollary}

\noindent
\textbf{Remark}: The notation ``$\rho$" is utilized in (\ref{def1}),  (\ref{eq_quadratic}) and (\ref{eq_quadratic1}). Unless explicitly stated otherwise, for the rest of the paper $\rho$ shall refer exclusively to the solution of (\ref{eq_quadratic1}). The primary QCQP (\ref{eq_quadratic1}) plays a crucial role in this paper. Not only does its solution, $\rho$, measure the observability of individual variables, but the concept and properties of the effective region/subspace introduced in Section \ref{sec_obs_effective} also rely on the numerical solution of (\ref{eq_quadratic1}). 




\subsection{Nonlinear systems}
\label{sec_nonlinear}
Corollary \ref{corollary1} is applicable to linear systems, which is not the case in most real world applications. For nonlinear systems, we formulate the primary QCQP (\ref{eq_quadratic1}) based on an empirical observability Gramian, an approximation of $G$ for nonlinear systems \cite{kreneride2009,kangxu2011}. Consider a nonlinear system (\ref{eq_sys}). Let $\{ \bfu (k) \}_{k=0}^K$ be a solution. The empirical observability Gramian is defined as follows. Let $h>0$ be a small number, let $\{ \hat \bfu_i(k) \}_{k=0}^K$ be the solution of (\ref{eq_sys}) with initial state 
\EQ
\hat \bfu_i(0)=\bfu (0)+h \bfe_i ,\nonumber
\EE 
where $\bfe_i$ is the $i$th unit vector in $\Real^n$. We define the empirical observability Gramian, $G$, as follows
\EQ
\label{eq_empiricalG}
&D_i\bfy(k) =  (\bfh(\hat \bfu_i(k))- \bfh(\bfu(k)))/h, & i\in \calI ,\\
&D\bfy (k)= \MT D_{i_1}\bfy(k)&D_{i_2}\bfy(k)&\cdots &D_{i_{n_\calI}}\bfy(k)\EM,\\
&G=\ds\sum_{k=0}^K \left( D\bfy(k)\right)^\intercal D\bfy(k),
\EE
where $\calI$ is a subset of $\{1, 2, \cdots,n\}$ and $n_\calI$ is the number of elements in $\calI$. The empirical observability Gramian is a $n_\calI\times n_\calI$ matrix. If $n_\calI =n$, then $G$ has full size, $n\times n$. For high dimensional problems, we choose an index subset so that $n_\calI < n$ to reduce the computational load. We will introduce in Section \ref{sec_obs_effective} a method of determining the index set $\calI$ that can significantly reduce the dimension of $G$. The matrix $W(K)$ in (\ref{eq_quadratic1}) is computed using finite difference,
\EQ
&D_iw =  (w (\bfu_i(K))- w(\bfu(K)))/h,\\
&W (K)= \MT D_{i_1}w&D_{i_2}w&\cdots &D_{i_{n_{\calI}}}w \EM, & \;\; i_j\in \calI .
\EE
Based on the empirical observability Gramian defined in (\ref{eq_empiricalG}), we can compute $G_\delta$ using (\ref{eq_Gmodified}). Then the primary QCQP (\ref{eq_quadratic1}) can be formulated and solved. The resulting $\rho$ is a quantitative measure of the observability of $z(K)$. 

\section{Effective region}
\label{sec_obs_effective}
For high-dimensional problems, such as the atmospheric models in numerical weather prediction, the computational burden of finding all eigenvalues of the empirical observability Gramian for a global model is beyond the capabilities of current and projected computers. Moreover, for limited-area estimation, obtaining the required lateral boundary conditions for LAMs typically involves using a lower-resolution global model. This approach increases the computational burden. In addition, lateral boundary conditions are recognized as a source of model error \cite{giorgimearns1999}. In this section, we introduce the concept of effective region/subspace. This concept serves two goals. Firstly, the computation of observability can be restricted to the effective subspace. As a result,  the dimension of the empirical observability Gramian is significantly lower than that of the full state space. Therefore, the concept of observability is computationally feasible for high dimensional systems. Secondly, we will demonstrate in Section \ref{sec_learning} that the surrogate DA model is independent of lateral conditions around the area of interest, provided that the area of interest is inside an effective region. This advantageous property makes the surrogate DA model fundamentally different from conventional data assimilation methods, such as EnKF and 4D-Var, which rely on the lateral boundary condition of LAMs in the estimation process. 

Throughout the remainder of this paper, we will concentrate on dynamical systems that result from the discretization of PDEs. Estimating the state and parameters for this kind of systems is important in a wide spectrum of applications, such as numerical weather prediction, climate change, and combustion dynamics. If the size of the grid is large, the computational burden of evaluating $G_\delta$ is formidable. A remedy that we propose in (\ref{eq_empiricalG}) is to perform the computation for $i$ in a small subset $\calI\subseteq \{1,2,\cdots,n\}$. Because the dimension of $G$ equals $|\calI|$, it is important to select a small $|\calI|$ so that the computational load can be significantly reduced, while the Gramian is still effective in determining the observability. It will be shown in this section that we can find $\calI$ based on the concept of the effective region. 

\begin{figure}[!ht]
\centering
\includegraphics[width = 3.5in]{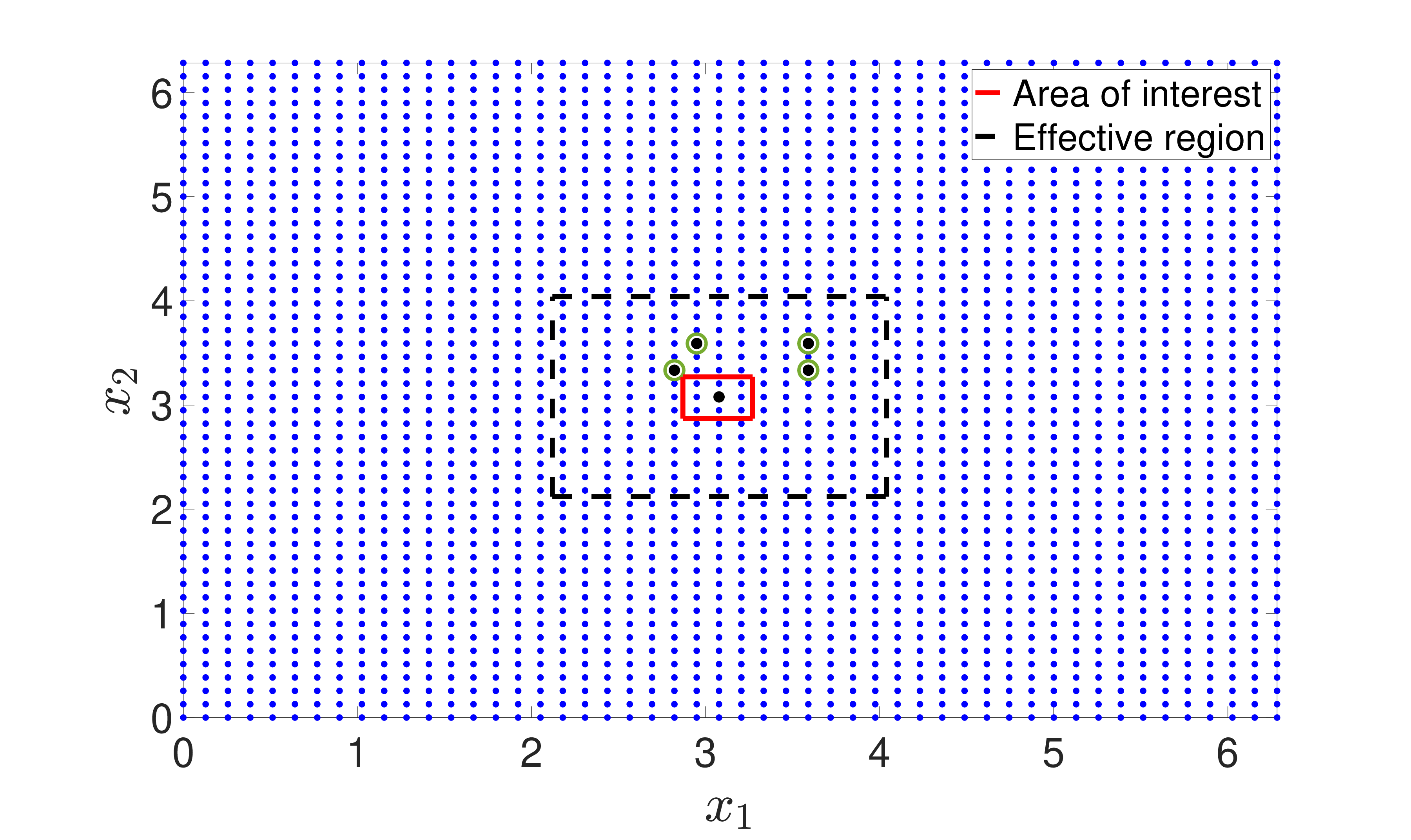}
\begin{minipage}{4.5in}
\caption{ An example of a 2D space. The sensors are located at the large dots. The variables to be estimated are the state at the grid points inside the red square. }
\label{fig_grid2}
\end{minipage}
\end{figure}

For illustrative purposes, let's consider the following Burgers equation in $\Real^2$, 
\EQ
\label{eq_burgers}
&\Fr{\partial U}{\partial t}+U\Fr{\partial U}{\partial x_1}+V\Fr{\partial U}{\partial x_2}=\kappa \left( \Fr{\partial^2 U}{\partial x_1^2}+\Fr{\partial^2 U}{\partial x_2^2} \right),\\
&\Fr{\partial V}{\partial t}+U\Fr{\partial V}{\partial x_1}+V\Fr{\partial V}{\partial x_2}=\kappa \left( \Fr{\partial^2 V}{\partial x_1^2}+\Fr{\partial^2 V}{\partial x_2^2} \right),
\EE
where $(x_1, x_2)\in [0, 2\pi]\times [ 0, 2\pi]$, $t\in [0, T]$. To numerically solve the PDE, we discretize the equation over a grid. Notation-wise, the grid point in space with coordinates $(i\varDelta x_1, j\varDelta x_2)$ is denoted by $(i,j)$.  The value of $U(x_1,x_2,t)$ at a grid point $(i,j)$ at time $t=k\varDelta t$ is denoted by $u_{i,j}(k)$. Similarly, we define $v_{i,j}(k)$. The discretized approximation of (\ref{eq_burgers}) is a system that can be represented in the form of (\ref{eq_sys}). The state variable is a vector in $\Real^n$, where $n=2(n_{x_1}\times n_{x_2})$, and $n_{x_1}$ and $n_{x_2}$ are the numbers of grid points in the $x_1$ and $x_2$ directions, respectively. Figure \ref{fig_grid2} shows a typical grid in which the large dots represent sensor locations, and the small box in the center is the area of interest. In general, it is not necessary to assume that the sensors are located at grid points. The concept in this section is applicable to any type of sensors that have an observation operator in the form of (\ref{eq_sysb}). The purpose of a limited-area DA algorithm is to estimate $u_{i,j}$ and $v_{i,j}$ at the grid points located inside the area of interest.

To evaluate the observability of $u_{i,j}(K)$ at a given time $K$, we solve the primary QCQP (\ref{eq_quadratic1}) and compute $\rho$, which has a closed-form solution (\ref{eq_tilderho}). However, the challenging issue is that computing the empirical observability Gramian (\ref{eq_empiricalG}) is computationally too expensive if the index set $\calI $ is too large. In the following, we introduce a way to numerically find $\calI$. 
In this example, $\kappa=0.14$ and $T=5$. The initial condition is
\EQ
\label{eq_initial_nominal}
&U(x_1,x_2,0)=V(x_1,x_2,0)=g(x_1,x_2),\\
&g(x_1,x_2)=\left\{ \begin{array}{lll} x_1^3 (2-x_1)^3 x_2^3 (2-x_2)^3, & (x_1, x_2)\in [0, 2]\times [0, 2]\\ 0 &\mbox{otherwise} \end{array}\right.
\EE
The equation is discretized on a uniform grid that has $50$ points in both $x_1$ and $x_2$ directions. Numerical scheme exploited is a first order upwind in time and space for the first order partial derivatives on the left-hand side of the equation, with an implicit central difference for the second order spatial derivatives on the right-hand side \cite{leveque}. The time step size is $\varDelta t=0.05$. The sensor error bound is set to be $\epsilon = 0.065$, which is about $8\%$ of the variation of sampling trajectories. The error upper bound of initial guess is $\delta=0.5$, which account for about $60\%$ of the variation of sampling trajectories. The step size in the computation of the empirical observability Gramian is $h=10^{-5}$. The parameters are summarized in Table \ref{table1}.
\begin{table}[!ht]
\begin{center}
\begin{tabular}{ |c|c|c|c|c|c| }
\hline
  $\kappa$ & $\epsilon$ & $\delta$ & grid & $\varDelta t$ & $h$ \\
  \hline
  0.14 & 0.065 & 0.5 &50&0.05&$10^{-5}$\\
  \hline
\end{tabular}
\caption{Parameters used in all examples in Section \ref{sec_learning_example}.}
\label{table1}
\end{center}
\end{table}

\begin{figure}[!ht]
\centering
\includegraphics[width = 3.5in]{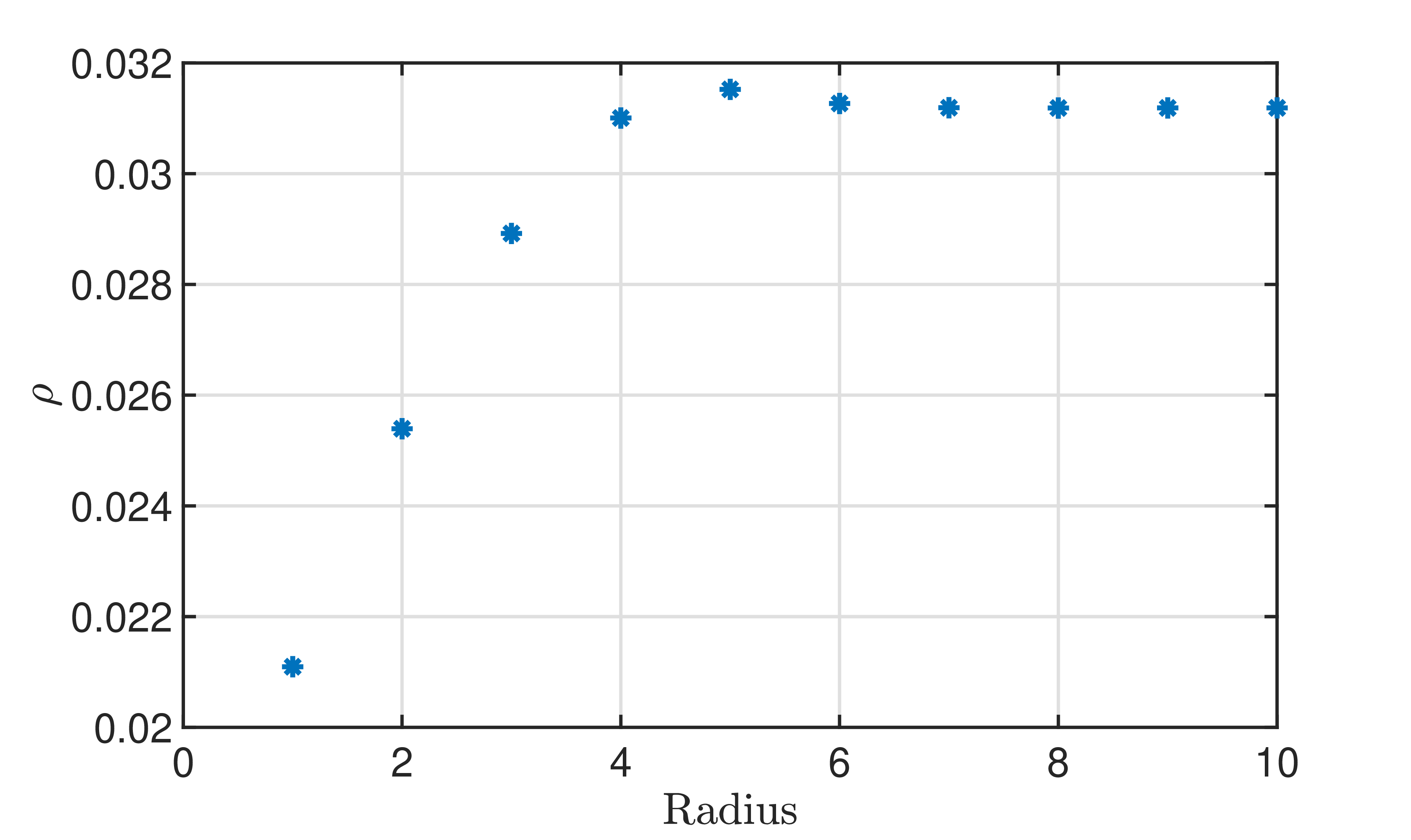}
\begin{minipage}{4.5in}
\caption{ The value of $\rho$ computed using different radii. }
\label{fig_obs_radius}
\end{minipage}
\end{figure}

The grid point located in the center of the small square in Figure \ref{fig_grid2} is $(25,25)$. For illustrative purposes, in the following we compute the value of $ \rho$ for $u_{25,25}$ along the nominal solution with the initial condition (\ref{eq_initial_nominal}). Assume that four sensors are used and located at
\EQ
(24,29),(29,27), (29,29),(23,27), 
\EE
as depicted in Figure \ref{fig_grid2}. The time window is $K=12$. 
Referring to (\ref{eq_empiricalG}), the empirical observability Gramian is computed based on the variations of a subset of state variables with indices in the set $\calI$. In this example, $\calI$ consists of all integer pairs, $(i,j)$, that represent the grid points located inside a square region around the center point $(25,25)$. The radius of the region is the number of grid points from $(25,25)$ to the boundary. For example, the radius of the square bounded by the dotted line in Figure \ref{fig_grid2} is R = $7$. The value of $ \rho$ is computed through solving the primary QCQP (\ref{eq_quadratic1}) in a sequence of regions with various radii. The result is shown in Figure \ref{fig_obs_radius}. The value of $ \rho$ is stabilized when R $\geq 7$. This implies that the observability of $u_{25,25}$ is determined by the empirical observability Gramian inside a region that is much smaller than the original spatial domain. 

\begin{definition}
Consider a system in the form of (\ref{eq_sys}). If the value of $ \rho$ computed based on $\calI$ in (\ref{eq_empiricalG}) is stabilized, i.e., the value of $ \rho$ computed based on any larger index set that includes $\calI$ has negligible change, then the state variables with indices in $\calI$ form a subspace of the state space. This subspace is called an effective subspace. For PDEs, if the grid points associated with the state variables in the effective subspace form a region, then it is called an effective region. 
\end{definition}

By this definition, the dotted square in Figure \ref{fig_grid2}  with radius R = $7$ is an effective region. If the effective region is much smaller than the domain of the PDE, the dimension of the empirical observability Gramian is significantly reduced. This is critical for high dimensional problems, such as the atmospheric models in numerical weather prediction, for which the eigenvalues of $G$ in its full dimension is impossible to compute. However, we would like to point out that finding an effective region is not a deterministic process. In this example, we search for the effective region by evaluating a sequence of $ \rho$ in a series of squares that increases in size until the value of $\rho$ does not change much. In this example, we consider the value of $\rho$ being stabilized if its variation is less than $10^{-4}$. This ad hoc approach is based on the assumption that the dynamics around two grid points that are far away do not affect each other within a limited time window. One may choose a larger effective region with R $> 7$. In some cases, as exemplified in the next section, a larger effective region improves DA accuracy, albeit at the cost of increased computational burden during the generation of data for training neural networks.

\section{Supervised learning for limited-area DA}
\label{sec_learning_example}
In this section, we introduce a surrogate DA model. It is exemplified by a discretized Burgers equation. 
\label{sec_learning}
\subsection{The architecture of surrogate DA model}
\label{sec_DAmodel}
We use the Burgers equation (\ref{eq_burgers}) to illustrate the design of a surrogate DA model. Similar to the previous section, the PDE is discretized on a uniform grid that has $50$ points in both $x_1$ and $x_2$ directions. For reasons of simplicity,  we assume that sensors are located at grid points. The measured data at a grid point $(i,j)$ forms a vector 
\EQ
\label{eq_Xij}
X_{i,j}=\left[ u_{i,j} (k_0) , v_{i,j} (k_0), \cdots, u_{i,j} (k_0+K), v_{i,j} (k_0+K) \right]^\intercal+{\boldsymbol\nu}_{i,j} ,
\EE
where $\boldsymbol\nu_{i,j}$ represents the censor uncertainty, which is a vector with an appropriate dimension. Its elements are random variables following an i.i.d. Gaussian distribution with zero mean and standard deviation $\sigma=0.065$. The initial time of the data point is $t_0=k_0\varDelta t$. 

The goal of the limited-area DA is to estimate the value of the state at the grid points located inside the red square in Figure \ref{fig_grid2}. In general, when designing a surrogate DA model, one can choose a subset of variables for estimation.  It is not necessary to include all state variables of the system when calculating the DA model. For illustrative purposes, we focus on the estimation of $U(x_1, x_2, t)$, i.e., 
\EQ
\label{eq_uij}
\{ u_{i,j}(k_0+K); 24\leq i, j\leq 26\}.
\EE
The process of estimating $V(x_1,x_2,t)$ is similar.  A data point for training and validation is a pair $(X,Y)$, where  
$X\in \Real^{2*N_{sensor}(K+1)}$ is formed by stacking  $X_{i,j} $ at all sensor locations, $N_{sensor}$ is the number of sensors, $Y\in \Real^9 $ consists of the nine variables in (\ref{eq_uij}). We first computed  $10,000$ solutions of the discretized Burgers equation in $t\in [0, T]$, $T=5$, with a constant zero boundary condition and the following random initial values around the nominal solution,
\EQ
\label{eq_random_initial}
&U(x_1,x_2,0)=g(x_1,x_2)+\ds\sum_{l=1}^{N_F}\sum_{s=1}^{N_F} a^U_{ij}\sin(\Fr{lx_1}{2})\sin(\Fr{sx_2}{2}),\\
&V(x_1,x_2,0)=g(x_1,x_2)+\ds\sum_{l=1}^{N_F}\sum_{s=1}^{N_F} a^V_{ij}\sin(\Fr{lx_1}{2})\sin(\Fr{sx_2}{2}),
\EE
where $a^U_{ij}$ and $a^V_{ij}$ have a normal distribution with a mean of zero and standard deviation of $0.1$, and we set $N_F=3$. 
The data set for training is generated from these solutions at uniformly distributed random time $t_0=k_0\varDelta t$ to form $(X,Y)$ in (\ref{eq_Xij}) and (\ref{eq_uij}).  A total of $30,000$ data points is generated.  Please note that the value of $k_0$ varies randomly for each data point. However, for simplicity, we consistently use the notation $k_0$ without distinction. The data set for validation is generated in the same way. It is worth noting that this is not an efficient way of generating data. In Section \ref{sec_learning_effectiveregion}, we introduce a more efficient method that utilizes a model with significantly reduced dimensions.

The surrogate DA model is a feedforward neural network that has $2N_{sensor}(K+1)$ inputs, the dimension of sensor data, and nine outputs. Given a sensor data point $X$, the output $Y$ of the neural network estimates $u_{i,j}(k_0+K)$ in (\ref{eq_uij}).
The neural network has $8$ hidden layers with a width of $16$ neurons. Each neuron is a hyperbolic tangent. In the supervised learning, the neural network is trained using the standard mean square error as the loss function. 

We numerically tested four different cases with varying numbers of sensors and different values of $K$.
\begin{table}[!ht]
\begin{center}
\begin{tabular}{ |c|c|c|c|c|c|c|} 
 \hline
$N_{sensor}$ & 8 & 4 & 4 &4  \\ 
 \hline
K & 12 & 12 &5 &2  \\ 
  \hline
  Average $\rho$ & 0.0209 & 0.0300 &0.0672 &0.1416  \\ 
  \hline
RMSE & 0.0095 & 0.0318 & 0.0445&0.0575  \\ 
  \hline
\end{tabular}
\caption{The RMSE of the surrogate DA model.}
\label{table2}
\end{center}
\end{table}
The estimation error of $u_{25,25}(k_0+K)$ is summarized in Table \ref{table2}. The third row shows the average value of $ \rho$ over $250$ different solutions with random initial conditions. Its value increases as the amount of sensor information decreases. In the case of eight sensors and $K=12$, the average $\rho$ is $0.0209$, which is about $2.6\%$ of the largest variation of $u_{25,25}$ in the data. This case is considered strongly observable. On the other hand, the average $ \rho =0.1416$ when $N_{sensor}=4$ and $K=2$ is about $17\%$ of the largest variation of $u_{25,25}$. This is considered practically unobservable. In Table \ref{table2}, the root-mean-square error (RMSE) of the surrogate DA model, which is the trained neural network, is consistent with the observability in the sense that the estimation error escalates as the value of $\rho$ increases. 

It is worthy to point out that estimators such as Kalman filter or particle filters take into consideration of estimation error's probability distribution, in addition to the sensor data. It is expected that these filters should achieve an error that is smaller than $ \rho$, which is defined based on the dynamical model but independent of the propagation of the uncertainty's probability distribution. On the other hand, an estimator can be deemed inadequate if its error surpasses $ \rho$. This is because $ \rho$ represents the worst estimation error among indistinguishable estimates.

\subsection{Learning inside the effective region}
\label{sec_learning_effectiveregion}
In Section \ref{sec_DAmodel}, the data sets are generated by solving the discretized PDE across the entire domain with a zero boundary condition.  However, this approach raises two important considerations. Firstly, is it feasible to generate data sets by solving the PDE within a limited area, rather than the entire domain, in order to significantly reduce the computational demand? Secondly, how well does the surrogate DA model generalize to different boundary conditions? For instance, if a surrogate DA model is trained using data sets with a zero boundary condition, how effectively does it perform when applied to data sets with non-zero boundary conditions? 

An intriguing characteristic  of the surrogate DA model is its insensitivity to boundary conditions.  By the definition of the effective region, $\rho$ can be determined by computing the empirical observability Gramian solely within the effective region, without the necessity of exploring state variations across the entire domain. Notably, alterations in the boundary conditions around the effective region do not alter the worst indistinguishable estimates. This property holds great significance as it enables us to generate data by utilizing a fabricated boundary condition around the effective region. It eliminates the requirement of providing precise lateral conditions, as commonly needed in LAMs utilized in applications such as numerical weather prediction.

For the purpose of illustration, let's consider the effective region computed for $u_{25.25}$. The radius of the region is R = $7$ (see Figure \ref{fig_grid2}), at which the value of $\rho$ is stabilized. The discretized Burgers equation is solved in the effective region, with random initial condition defined in (\ref{eq_random_initial}). Without solving the PDE across the entire domain, we do not have accurate boundary condition around the effective region. For computation, we introduce a fabricated boundary condition,
\EQ
\label{eq_fabricate_IC}
U(x_1,x_2,t)=U(x_1,x_2,0)e^{-t},\\
V(x_1,x_2,t)=V(x_1,x_2,0)e^{-t}.
\EE
The data sets are generated following the same method described in Section \ref{sec_DAmodel}, with the only difference being that we solve the Burgers equation within the effective region using the fabricated boundary condition in (\ref{eq_fabricate_IC}). The effective region of $u_{25,25}$ consists of a grid with dimensions $15\times 15$. Any grid points outside this region, regardless of their quantity, are irrelevant. By concentrating on the effective region, the computational requirements can be significantly reduced by eliminating the need to solve the partial differential equation in a substantially larger spatial domain. In addition, the trained surrogate DA model is insensitive to boundary conditions around the effective region, i.e., it is applicable to solutions of the Burgers equation even if they do not satisfy the fabricated boundary condition (\ref{eq_fabricate_IC}). In the following, we use a variety of different boundary conditions to verify this property. 

\begin{figure}[!ht]
\begin{subfigure}{.5\textwidth}
\includegraphics[width = 3.3in]{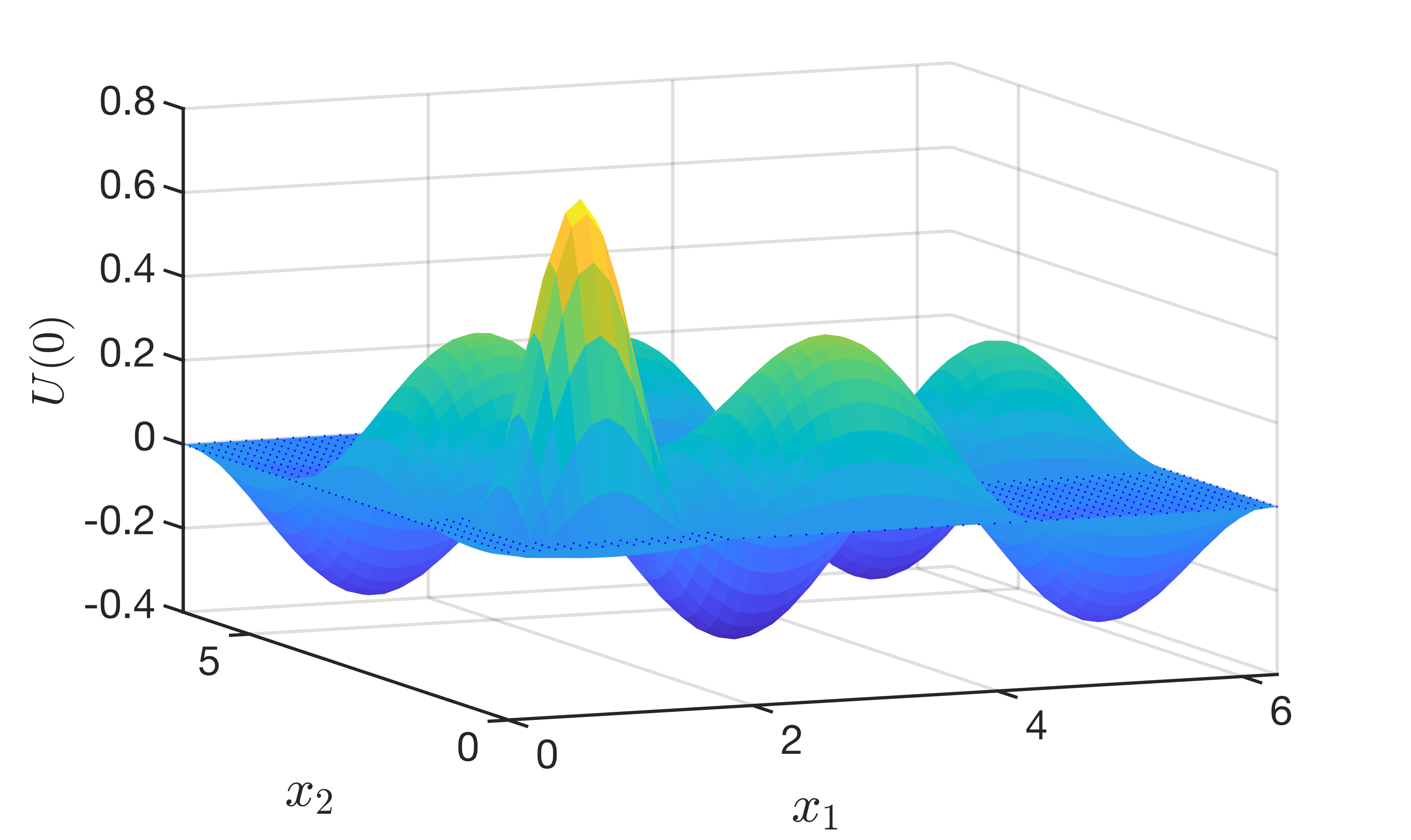}
\caption{}
\end{subfigure}
\begin{subfigure}{.5\textwidth}
\includegraphics[width = 3.3in]{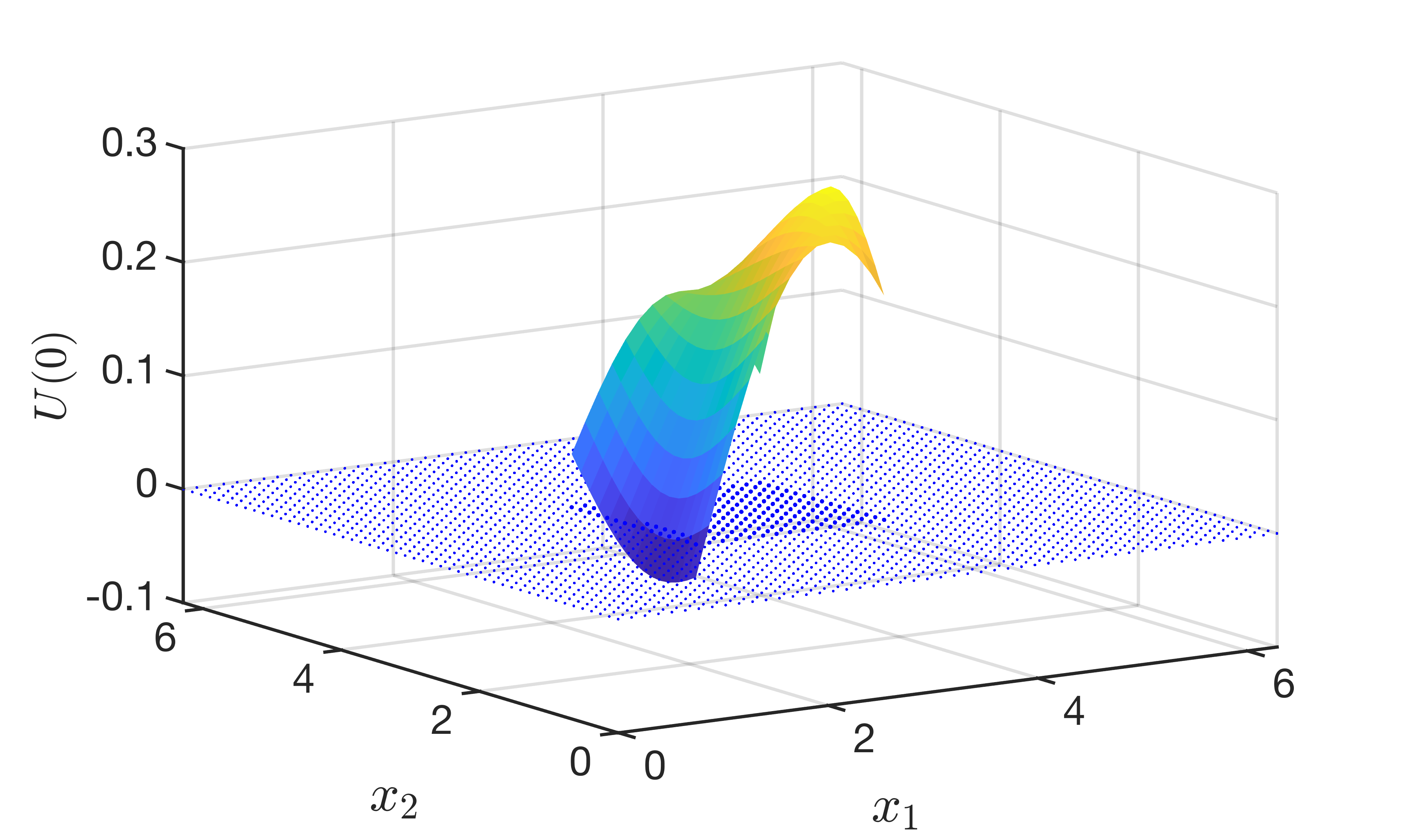}
\caption{}
\end{subfigure}
\caption{ An example of initial condition, validation data in (a) and training data in (b). }
\label{fig_IC}
\end{figure}

The validation data set generated in Section \ref{sec_DAmodel} has $30,000$ samples. Their initial conditions are defined over the original spatial domain of the PDE with $50\times 50$ grid points. A zero boundary condition is applied. This is significantly different from the training data that is generated in the effective region  with the fabricated boundary condition. Two examples of typical initial conditions are shown in Figure \ref{fig_IC}. The differences are obvious. We use the training data from the effective region to train surrogate DA models for different number of sensors with different length of time windows. Then the surrogate DA model is applied to the validation data that does not follow the fabricated boundary condition. We assume that the sensor uncertainty is a random variable following an i.i.d. Gaussian distribution with zero mean and standard deviation $\sigma=0.065$. The RMSE of the surrogate DA model is summarized in Table \ref{table3}

In comparison to Table \ref{table2}, the surrogate DA model trained with the fabricated boundary condition exhibits a slightly higher RMSE. Nonetheless, the overall estimation is considered satisfactory as the RMSE values are either lower than or comparable to the value of $ \rho$. Moreover, the RMSEs are smaller than the standard deviation of sensor uncertainty, $\sigma = 0.65$, indicating that the surrogate DA model effectively filters out the sensor noise.
\begin{table}[!ht]
\begin{center}
\begin{tabular}{ |c|c|c|c|c|c|c|} 
 \hline
$N_{sensor}$ & 8 & 4 & 4 &4  \\ 
 \hline
K & 12 & 12 &5 &2  \\ 
  \hline
  Average $\rho$ & 0.0209 & 0.0300 &0.0672 &0.1416  \\ 
  \hline
  RMSE (R = $7$) & 0.0129 & 0.0425 & 0.0421&0.0497  \\
  \hline
\end{tabular}
\begin{minipage}{4.5in}
\vspace{0.1in}

\caption{The RMSE of the surrogate DA model, trained using data exclusively from the effective region R = $7$, but validated using solutions solved in the original domain.}
\label{table3}
\end{minipage}
\end{center}
\end{table}
\vspace{-0.3in}

If we use a radius at which $\rho$ has not yet stabilized, for instant R = $4$, the variation of the boundary condition around this region has a greater impact on the solution of the PDE compared to the case when R = $7$. Consequently, the estimation error of the surrogate DA model, trained using the fabricated boundary condition, can be higher. We conducted a numerical test specifically for R = $4$.   The RMSE of the surrogate DA model  is summarized in Table \ref{table3b}.  When comparing these results to those in Table \ref{table3}, the error increases by $24\%$ if $N_{sensor}=8$ and $K=12$. The increase reaches $93\%$ when $N_{sensor}=4$ and $K=5$.

\begin{table}[!ht]
\begin{center}
\begin{tabular}{ |c|c|c|c|c|c|c|} 
 \hline
$N_{sensor}$ & 8 & 4 & 4 &4  \\ 
 \hline
K & 12 & 12 &5 &2  \\ 
\hline
 RMSE (R = $4$) & 0.0160 & 0.0549 & 0.0813&0.0546\\
  \hline
\end{tabular}
\begin{minipage}{4.5in}
\vspace{0.1in}

\caption{The RMSE of the surrogate DA model, trained using data exclusively from the region with R = $4$, but validated using solutions solved in the original domain.}
\label{table3b}
\end{minipage}
\end{center}
\end{table}

\subsection{Non-zero boundary condition}
\begin{figure}[!ht]
\begin{subfigure}{.5\textwidth}
\includegraphics[width = 3.3in]{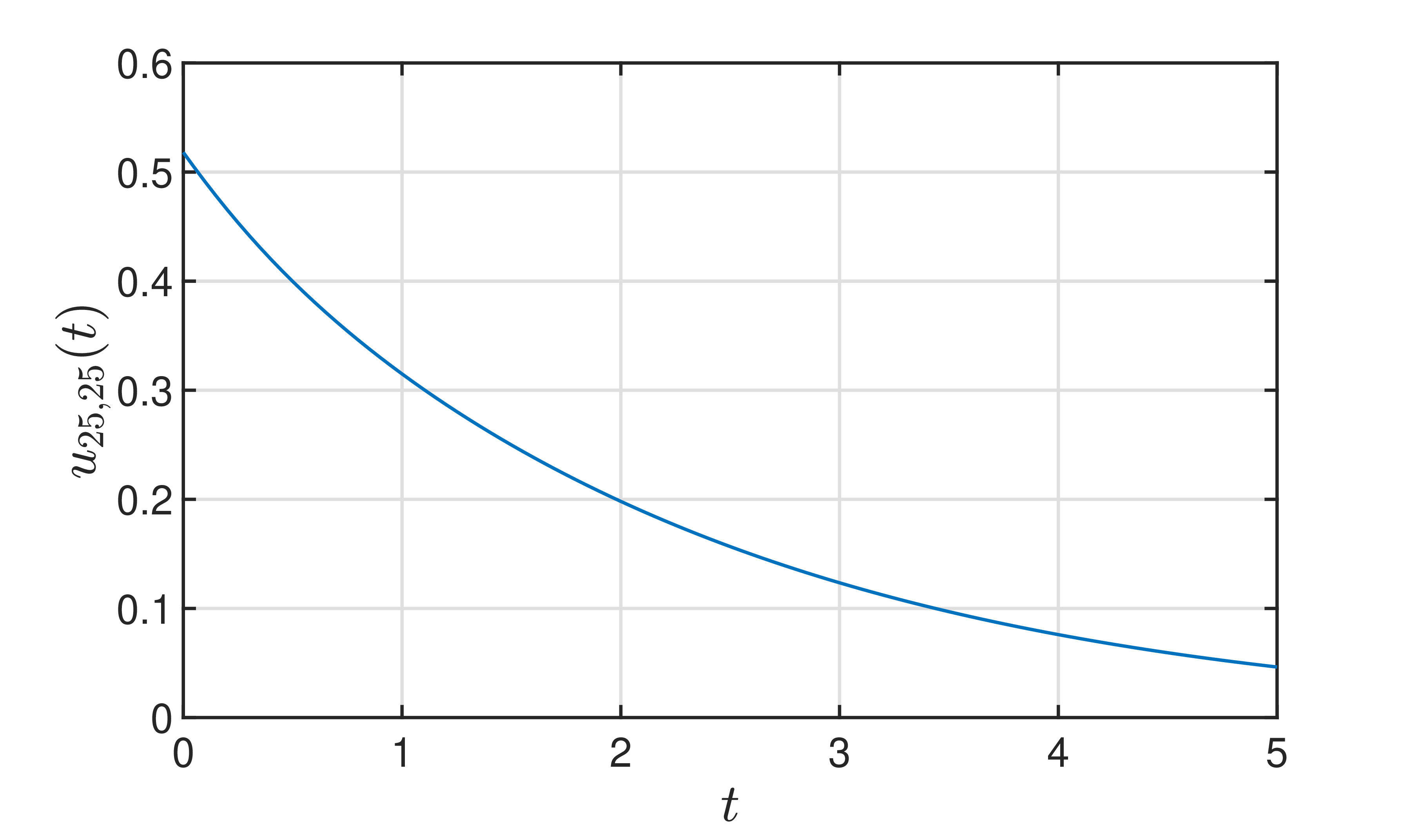}
\caption{}
\end{subfigure}
\begin{subfigure}{.5\textwidth}
\includegraphics[width = 3.3in]{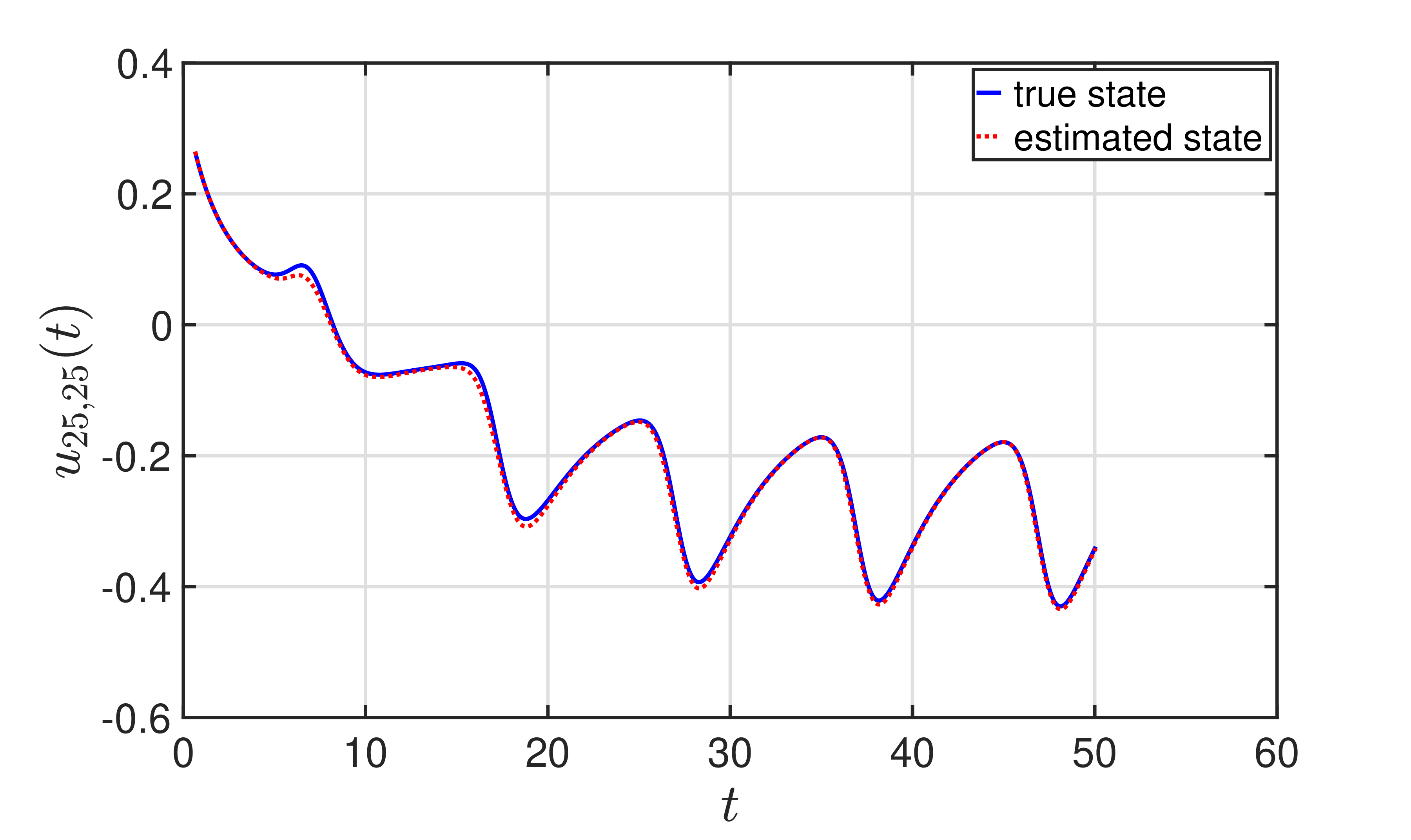}
\caption{}
\end{subfigure}
\begin{center}
\begin{minipage}{6.0in}
\vspace{-0.2in}

\caption{ (a) A typical trajectory of $u_{25,25}(t)$ computed in the effective region using a fabricated boundary condition. (b) A trajectory computed in the original domain with the boundary condition (\ref{eq_nonzerosboundary}). }
\label{fig_nonzeroboundary}
\end{minipage}
\end{center}
\end{figure}

\begin{figure}[!ht]
\centering
\includegraphics[width = 3.5in]{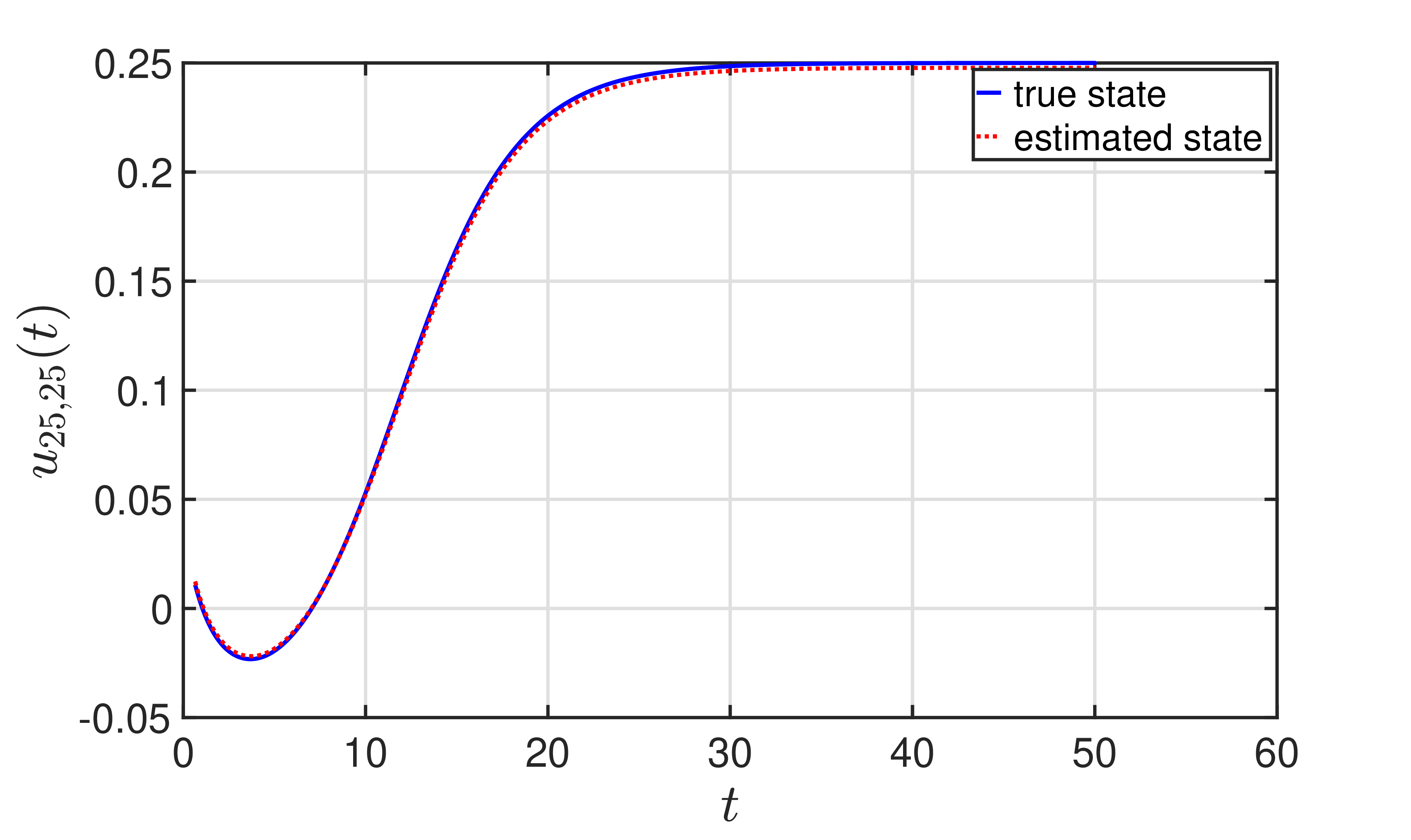}
\begin{minipage}{4.5in}
\caption{ The estimation of a trajectory with a constant boundary condition $u_{i,j}(t)=0.25\tanh (t)$.}
\label{fig_nonzerobundary2}
\end{minipage}
\end{figure}

In Section \ref{sec_learning_effectiveregion}, the surrogate DA model trained on data from the effective region with a fabricated boundary condition is validated using data from solutions solved across the entire domain. The validation data, however, equal zero along the boundary of the original spatial domain.  In this section,  we aim to investigate how well the surrogate DA models from Section \ref{sec_learning_effectiveregion} generalize to nonzero boundary conditions. The plot in Figure \ref{fig_nonzeroboundary}(a) is a typical trajectory of $u_{25,25}(t)$ computed in the effective region using the fabricated boundary condition. The surrogate DA model is trained using data from trajectories of this nature. We then apply the trained surrogate DA model to a trajectory shown in Figure \ref{fig_nonzeroboundary}(b). This trajectory is computed in the original spatial domain with the following nonzero boundary condition, 
\EQ
\label{eq_nonzerosboundary}
&U(z, t)=V(z,t)=\sin (\Fr{t\pi}{5})\sin(z),\\
&z=\left\{ \begin{array}{lll} x_1, & \mbox{ along horrizental boundaries} ,\\ x_2, & \mbox{ along vertical boundaries}.\end{array}\right.
\EE
Despite the visually dissimilar plots in Figure \ref{fig_nonzeroboundary}(a) and (b), the surrogate DA model generalizes well. The estimated trajectory is the dotted curve in Figure \ref{fig_nonzeroboundary}(b). The RMSE is $0.0068$ without sensor uncertainty and $0.0373$ if we add sensor uncertainty that follows an i.i.d. Gaussian distribution with standard deviation $\sigma=0.065$.  Similarly, the surrogate DA model performs well in the estimation of a trajectory with a constant non-zero boundary condition, as shown in Figure \ref{fig_nonzerobundary2}.  The RMSE is $0.0020$ without sensor uncertainty, and $0.0127$ otherwise. 

\begin{figure}[!ht]
\begin{subfigure}{.5\textwidth}
\includegraphics[width = 3.3in]{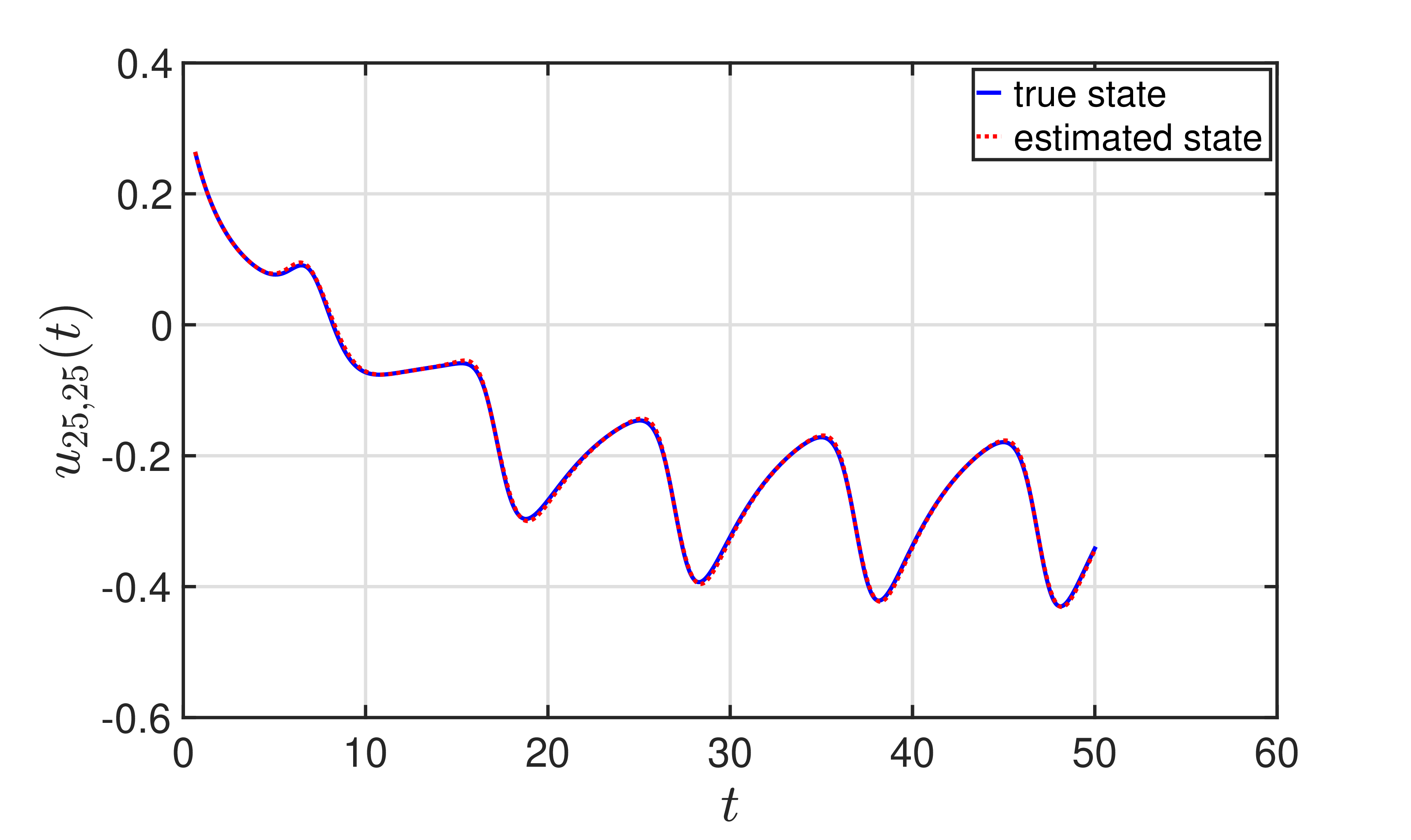}
\caption{}
\end{subfigure}
\begin{subfigure}{.5\textwidth}
\includegraphics[width = 3.3in]{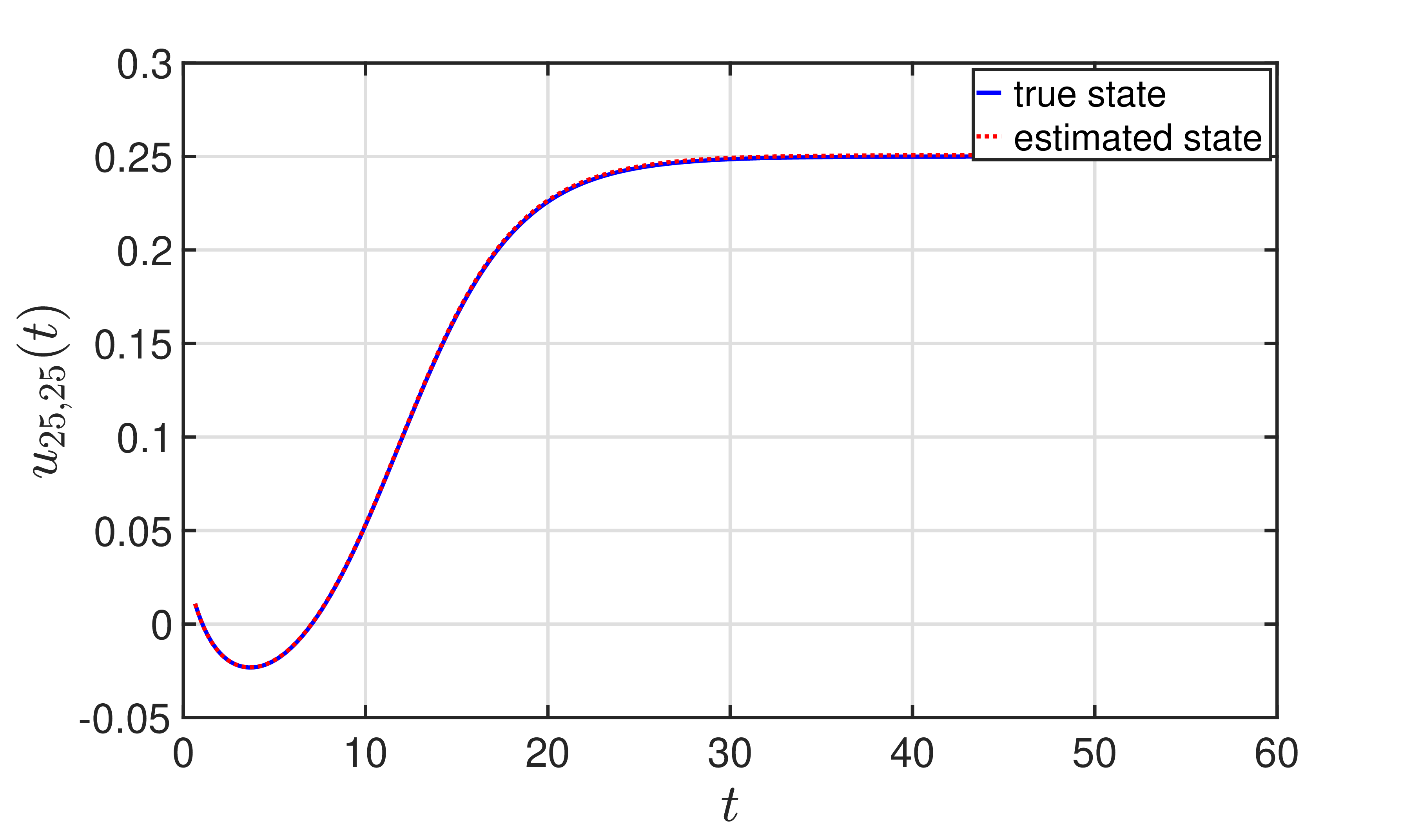}
\caption{}
\end{subfigure}
\begin{center}
\begin{minipage}{6.0in}
\vspace{-0.2in}

\caption{ Two trajectories that are estimated using the surrogate DA model trained in the effective region with R = $9$. In both (a) and (b), the estimated trajectory (dotted curve) follows the true trajectory (solid curve) closely.}
\label{fig_nonzeroboundary_R9}
\end{minipage}
\end{center}
\end{figure}

In this study, we discovered that enlarging the effective region leads to improved accuracy when generalizing the surrogate DA model to nonzero boundary conditions. For example, if we train the neural network in the effective region with R = $9$, the accuracy of the surrogate DA model is significantly improved. Specifically, the trajectory shown in Figure \ref{fig_nonzeroboundary_R9}(a) yielded an RMSE of $0.0027$ without sensor uncertainty and $0.0262$ otherwise, while the trajectory depicted in Figure \ref{fig_nonzeroboundary_R9}(b) exhibited an even lower RMSE of $6.4\times 10^{-4}$ without sensor uncertainty, and $0.0127$ with sensor uncertainty taken into account.

\begin{figure}[!ht]
\centering
\includegraphics[width = 3.5in]{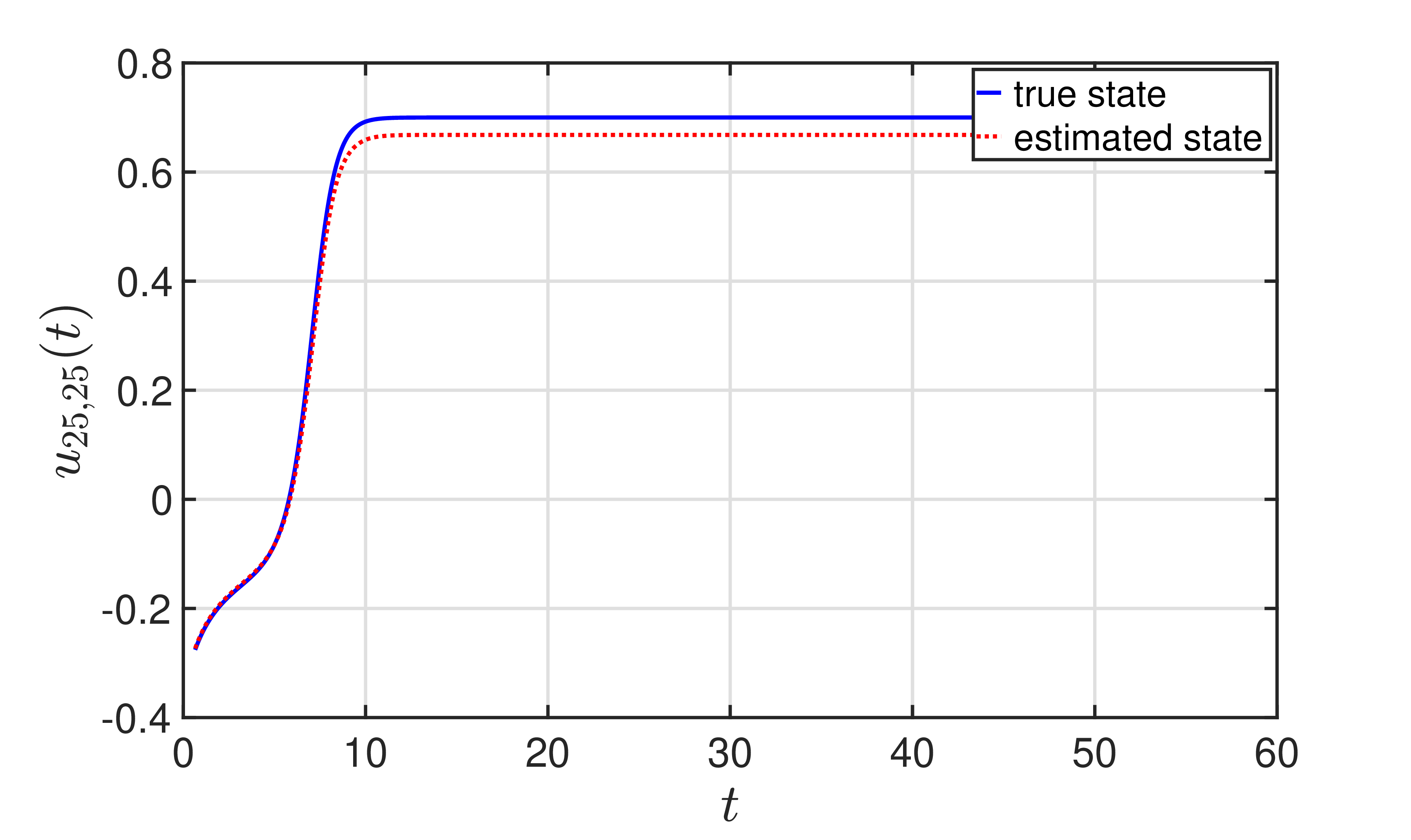}
\begin{minipage}{4.5in}
\caption{ The estimation of a trajectory with a constant boundary condition $u_{i,j}(t)=0.7\tanh (t)$.}
\label{fig_nonzerobundary3}
\end{minipage}
\end{figure}

We would like to emphasize that the surrogate DA model's insensitivity to boundary condition has its limitation. For instance, when the constant boundary condition is heightened, the surrogate DA model, which is trained using the fabricated boundary condition encompassing the effective region, may fail to accurately follow the true solution. An illustrative example can be observed in Figure \ref{fig_nonzerobundary3}. The surrogate DA model exhibits an underestimation of the stabilized value of $u_{25,25}(t)$. This phenomenon can be explained by the fact that the maximum value of $u_{25,25}(t)$ exceeds the range covered by the training data. In general, ensuring that the training data contains sufficient richness to represent the full dynamics of the system remains an unresolved aspect that warrants further research.

\subsection{Shifting the area of interest}
\begin{figure}[!ht]
\centering
\includegraphics[width = 3.5in]{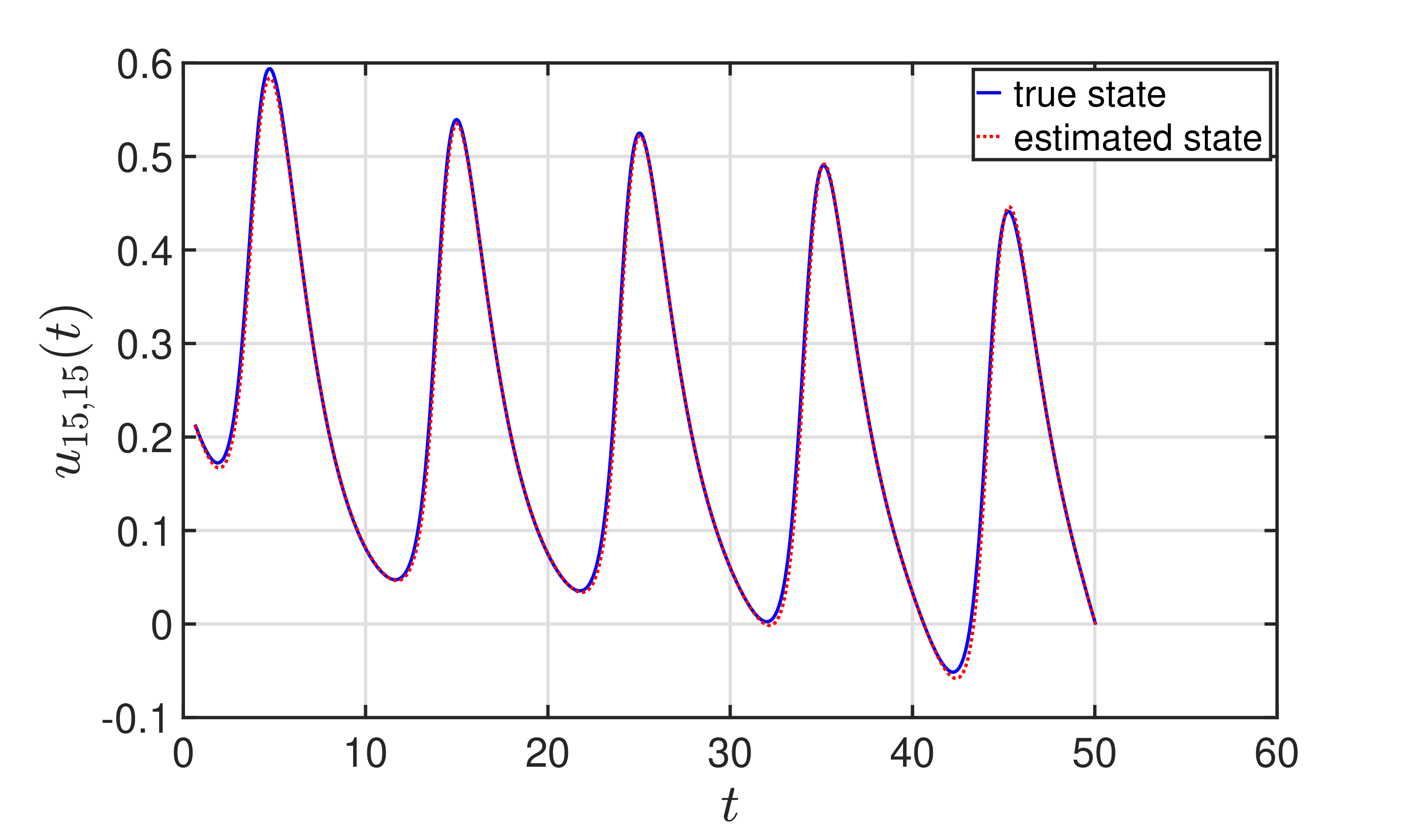}
\begin{minipage}{4.5in}
\caption{ The estimation of $u_{15,15}(t)$ by the surrogate DA model trained in the effective region around the grid point $(25,25)$ with R = $7$.}
\label{fig_nonzerobundary4}
\end{minipage}
\end{figure}

For PDEs in which coefficients are independent of the spacial variable, such as the Burgers equation (\ref{eq_burgers}), a surrogate DA model trained in one effective region can be applicable to another area as long as the relative positions of the point of interest and the sensor locations remain unchanged. To illustrate this, we consider the surrogate DA model presented in Section \ref{sec_learning_effectiveregion}, which is trained in the effective region depicted in Figure \ref{fig_grid2}. The center of the area of interest is located at the grid point $(25,25)$. Now, we employ the same surrogate DA model to estimate the state at a different location, $(15,15)$, assuming that the sensor locations have been shifted accordingly. Figure \ref{fig_nonzerobundary4} demonstrates the case of eight sensors with $K=12$. This trajectory is from a solution of the PDE that has a nonzero boundary condition in the original spatial domain.  Remarkably, despite being trained in a distinct effective region, the surrogate DA model performs well in estimating the trajectory $u_{15,15}(t)$, yielding an RMSE of $0.0092$. The RMSE is significantly reduced to $0.0027$ if we train the surrogate DA model in the effective region with R = $9$. 

\section{Learning-based prediction}
\label{sec_prediction}
\begin{figure}[!ht]
\centering
\includegraphics[width = 3.3in]{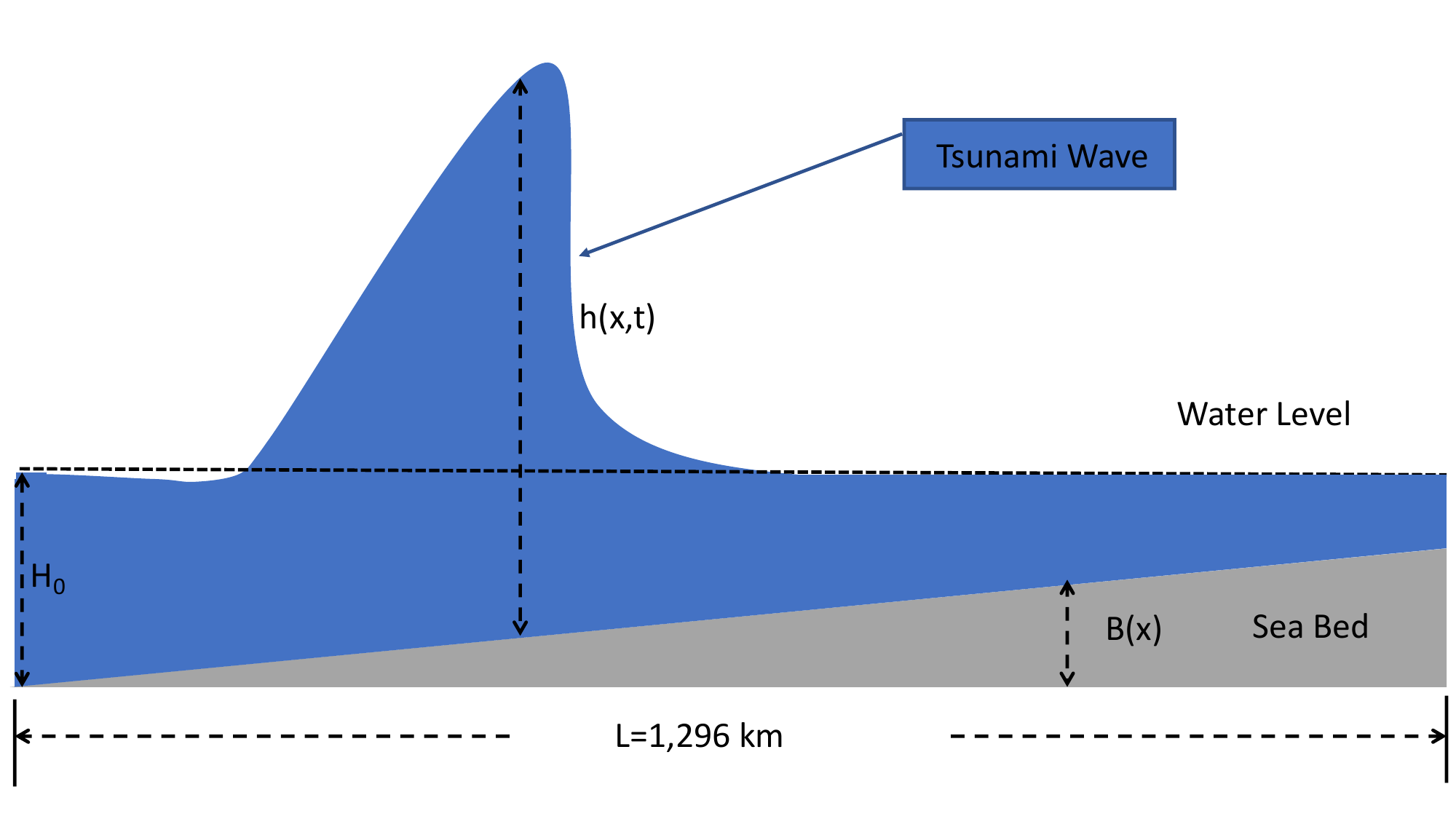}
\begin{minipage}{4.5in}
\caption{ The variables in the model of a tsunami wave. }
\label{fig_tsunami_variables}
\end{minipage}
\end{figure}

The same idea presented in Section \ref{sec_learning_example} is applicable to the prediction of dynamical systems. In this case, the prediction model is a neural network whose input is sensor information and the output is the predicted value of a state variable. In the following, we illustrate the idea using an example of tsunami forecast. 
The PDE dynamic model and parameters in this example  are adopted from  \cite{deleanu} in which a tsunami wave is simulated using the shallow water equations. Illustrated in Figure \ref{fig_tsunami_variables}, the water wave travels a horizontal distance of $L=1.296\times 10^{6}$ m. It approaches the shore on a variable seabed depth that has a constant slope $s=40/L$. The wave propagation is modeled using the shallow water equations
\EQ
&\Fr{\partial h}{\partial t}+\Fr{\partial \;}{\partial x}(uh)=0,\\
&\Fr{\partial \;}{\partial t}(uh)+\Fr{\partial \;}{\partial x}\left(v^2h+\Fr{1}{2}gh^2\right)=-gh\Fr{d B}{d x},
\EE
where $h=h(t,x)$ is the flow depth, $u=u(t,x)$ is the flow velocity, $B(x)$ is the seabed height, $g$ is the gravitational acceleration and $(t,x)\in [0, T]\times [0, L]$. In this example, the wave propagation is simulated for the time interval $T=60,000$ seconds. A total of $m=81$ sensors are  placed along the seabed with an equal distance of $16,200$m apart. Each sensor measures the full state $[h, uh]^T$. The value of parameters is summarized in Table \ref{table4}.

\begin{figure}[!ht]
\begin{center}
\includegraphics[width=3.3in]{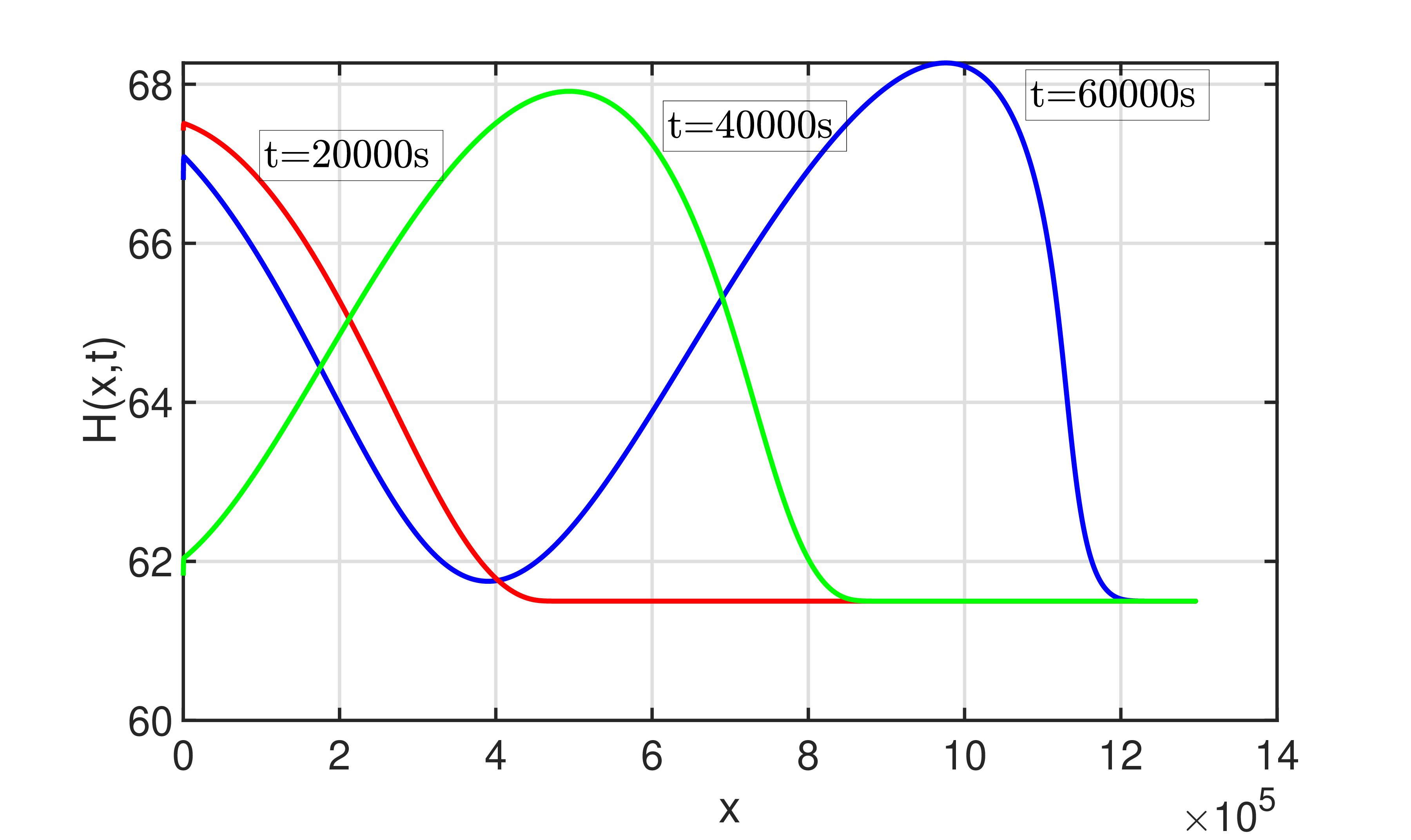}\includegraphics[width=3.3in]{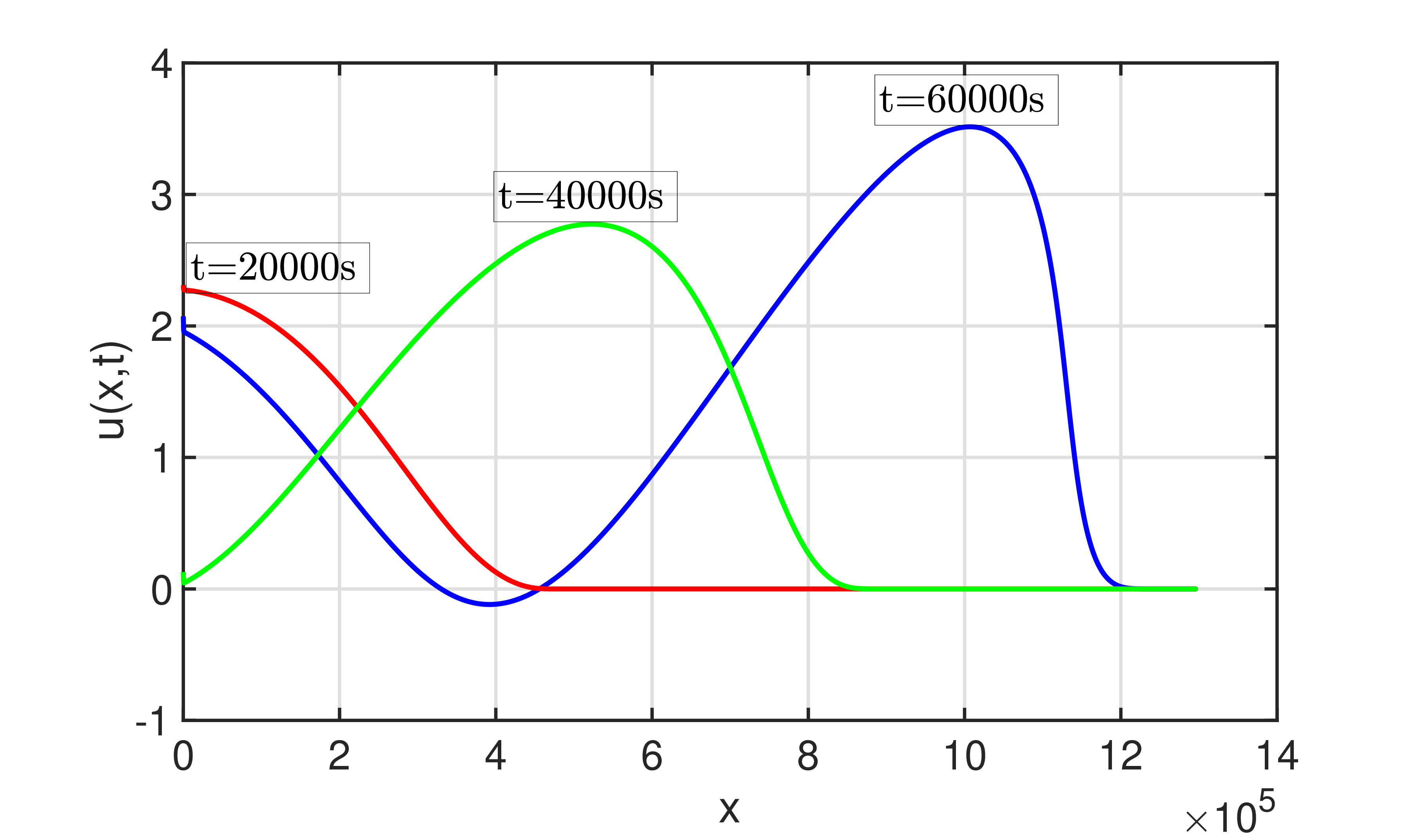}
\end{center}
\caption{ Water wave at $t=20000$, $40000$, $60000$ seconds.}
\label{fig_nominaltrajectory}
\end{figure}

\begin{table}
\begin{center}
\begin{tabular}{ |l |c| c| c| c|c| }
\hline
 Parameter & $L$ & $T$&$H_0$ & $M_h$ & g \\ 
 \hline
 Value & $1.296\times 10^6$m & $60,000$s & $61.5$m& $[0.2, 3.2]$ &$9.81$m/${\rm s}^2$   \\
 \hline
\end{tabular}
\caption{The parameters in the tsunami example.}
\label{table4}
\end{center}
\end{table}

The boundary conditions are given as follows,
\EQ
&B(x)=\Fr{40x}{L}, \\
&h(t,0)=64.5+M_h\sin\left(\pi\left( \Fr{4t}{86400} -\Fr{1}{2} \right) \right), && h(t,L)=H_0, \\
&u(t,0)=\sqrt{g}\left(\sqrt{h(t,0)}-\sqrt{H_0}\right), && u(t,L)=0, \\
&h(0,x)=H_0, &&u(0,x)=0.
\EE
Coded in Matlab, the shallow water equations are numerically solved using the first order Lax-Friedrichs algorithm \cite{leveque} over equally spaced grid points
$$x_0=0<x_1 < x_2 <\cdots < x_{N_x}=L, \;\; t_0=0<t_1<t_2<\cdots < t_{N_t}.$$
In the discretization, $\Delta x=324$m, $N_x=4,001$, $\Delta t=10$s and $N_t=6,000$. Shown in Figure \ref{fig_nominaltrajectory} is the nominal trajectory that we use to test a learning-based surrogate prediction model. 

The goal is to provide $900s$ prediction of the wave height at a high value facility located at $x=749,412$m. The training data is collected from a sequence of simulations in which $M_h$ takes $91$ equally space values in $[0.2,  3.2]$. A total of $9,100$ data points are randomly selected from the sample trajectories for neural network training and validation. Each data point is a pair $(X,Y)$ where $X$ contains the value of $h$ and $uh$ at all sensor locations to form a vector in $\Real^{162}$ at time $t$. Gaussian white noise is added to $X$, having a standard deviation of $0.15$. The scalar value of $Y$ equals $h$ at $t+9000$. The neural network has $162$ inputs and a single output. It has $8$ hidden layers with $16$ neurons in each layer. The surrogate prediction model performs well, yielding an RMSE of $0.0018$m and the maximum prediction error of $0.0283$m.

To measure the predictability of a variable, one can apply the same idea in observability analysis. More specifically, suppose we want to predict a variable at time $K+\varDelta K$ based on the observation $\{ \bfy (k)\}_{k=0}^K$. Replacing $z(K)$ and $W\bfu(K)$ in (\ref{def1}) or the primary QCQP (\ref{eq_quadratic1}) with the predicted values $z(K+\varDelta K)$ and $W\bfu(K+\varDelta K)$ yields a quadratic programming problem. The maximum value $\rho$ represents the worst error of indistinguishable predictions of $z(K+\varDelta K)$, thus serving as a measure of predictability. 

\section{Conclusions and discussions}
The proposed surrogate DA model for state estimation in a limited area does not require multiple integrations of a LAM during online computations. Instead, the online computation merely involves evaluating a feedforward neural network, which is significantly less computationally expensive compared to integrating high-dimensional LAMs. In the offline computation, which involves data generation and training of the neural network, the integration of the LAM does not require lateral conditions. The data is generated using a fabricated boundary condition. Furthermore, a reduced-order model can be employed to generate the dataset, focusing only on the effective region and avoiding the need for integrating the PDE across the entire domain. All these properties of the surrogate DA model offer clear advantages when compared to traditional limited-area DA methods. From a systems and control theory perspective, the notion of an effective region represents a novel application of observability. Moreover, we discovered a closed form solution to the primary QCQP in measuring the observability of linear systems. This formula can also be extended to nonlinear systems using the empirical observability Gramian. Importantly, the Gramian is computed specifically over the effective region, eliminating the need to solve the PDE across the entire spatial domain.

In this study, we acknowledge the presence of observation uncertainties in both the theoretical framework and illustrative examples. However, further investigation is required to understand the influence of model uncertainties on the effective region and surrogate DA model. Additionally, it is widely recognized that having information about the probability distributions associated with the uncertainties in both model and observation data is crucial for enhancing the accuracy of DA algorithms. Exploring the impact of uncertainty probability distributions on the effective region and neural network training in the context of limited-area DA presents numerous open questions that warrant future research. While the examples of the Burgers equation demonstrate encouraging outcomes, it is imperative to conduct further numerical experiments and explore diverse applications to validate the effectiveness and efficiency of the concepts and algorithms in real-life scenarios. Additional examples and applications are necessary to strengthen the evidence and ensure the practical applicability of the proposed methods.\\

\noindent
\textbf{Acknowledgement}\\
Any opinions, findings, and conclusions or recommendations expressed in this paper are those of the authors and do not necessarily reflect the views of the National Science Foundation and the Naval Research Laboratory. 

We express our sincere gratitude to Professor Anthony Austin for his invaluable discussions on topics related to computational linear algebra. These engaging conversations were inspirational in the derivation of the formula (\ref{eq_exp1}). 

\bibliography{myreference}
\bibliographystyle{ieeetr}







\end{document}